\documentclass{article}
\usepackage[a4paper,margin=2cm,centering,nohead,ignorefoot,footskip=5ex]{geometry}
\usepackage[english]{babel}
\usepackage{amsmath, amsxtra, latexsym,amscd, amsthm, amssymb,indentfirst}
\usepackage[mathscr]{eucal}
\usepackage{nccmath}
\usepackage{enumerate}
\usepackage{hyperref}
\usepackage{textcomp}
\usepackage{array, tabularx, longtable}%
\usepackage{multicol}
\usepackage{color}

\def\R{\mathbb R}
\def\N{\mathbb N}

\def\E{\mathbb E}

\def\int{{\rm int\,}}
\def\P{\mathcal P}

\def\P{\mathcal P}

\def\eps{\varepsilon}

\theoremstyle{plain}
\newtheorem{theorem}{Theorem}[section]
\theoremstyle{theorem}
\numberwithin{theorem}{section}

\newtheorem{proposition}[theorem]{Proposition}
\theoremstyle{definition}
\newtheorem{definition}[theorem]{Definition}

\newtheorem{example}[theorem]{Example}
\newtheorem{examples}[theorem]{Examples}
\newtheorem{remark}[theorem]{Remark}

\newtheorem{convention}[theorem]{Convention}
\begin{document}
\title{On non-linear optimization with a perturbed objective function}
\date{}
\maketitle
\centerline{\author{Nam Van Tran}}
\centerline{Faculty of Applied Sciences, Ho Chi Minh City University of Technology and Education,    Vietnam}
\centerline{namtv@hcmute.edu.vn}
\vspace{0,3cm}
\centerline{\author{Imme van den Berg}}
\centerline{CIMA, University of \'{E}vora,  Portugal} \centerline{ivdb@uevora.pt}
\vspace{0,5cm}
\begin{abstract}
\noindent A Lagrange multiplier theorem is derived for the case of an imprecise objective function and a precise constraint. The proof uses methods of analysis which deal  in a direct, algebraic way with imprecisions. They include imprecise differentiation, and an approximate Fermat Lemma and Implicit Function Theorem.
The tools are the external numbers of Nonstandard Analysis, which are models of Sorites imprecisions.
\vspace{0,5cm}

	\vspace{0,01cm}
	
\end{abstract}

\section{Introduction}

 The Lagrange multiplier method leads to necessary optimality conditions for non-linear optimization problems with constraints \cite{Karl Monty Toda}. In addition it enables sensitivity analysis, for the multiplier is an indicator for the effects of changes in the values of the constraint function \cite{Rockafellar}\cite{Sentivity}. In this article we extend the sensitivity analysis to the values of the objective function. We allow for small perturbations or errors in the form of imprecisions of Sorites-type, i.e. they are stable under small shifts, or some additions and multiplications. The principal result of this article (Theorem \ref{Lagrange multiplier}) concerns the existence of Lagrange multipliers for an optimization problem with an imprecise objective function and precise constraints. This relates to situations in which, due to all kinds of uncertain circumstances, the objectives are only known within some range of imprecisions, while the constraint, say a budgetary restriction, is rigid. In this context we show that an optimum can only be approximate, and indicate the size of the imprecisions occurring in the Lagrange equations for the multipliers. 
 
 Like in an earlier article on approximate linear programming \cite{Nam}, we model the imprecisions by the scalar neutrices and external numbers \cite{Koudjeti Van den Berg} of Nonstandard Analysis. These are (external) subsets of the reals, and permit an alternative for the Theory of Neglecting of Van der Corput \cite{Van der Corput} based on functional neutrices, which are generalizations of the  $ O(\cdot)'s $ and $ (\cdot)'s $. The calculation rules originate from informal Error Analysis \cite{Taylor} and come close to the rules for the real numbers, including a total order and a form of completeness. The resulting structure is called a Completely Arithmetical Solid \cite{dinisberg 2016}\cite{DinisBerg}. 
 
Functions from the reals to the external numbers bear some imprecisions and are called \emph{flexible functions}. We search for near-optimizers of flexible functions by approximate differentiation, and to this end we introduce new types of limits in terms of external numbers. This limits minimize the imprecisions of the sequential limits of \cite{Sequences} and permit to define forms of imprecise continuity and differentiation. We show that these notions satisfy many operations, including a Chain Rule. In the context of functions of two variables,  we define with these notions approximate partial derivatives, which lead to an approximate total differential. With the help of a Fermat Lemma for near-optimums and an Implicit Function Theorem we come to the Main Theorem on Lagrange multipliers for flexible functions, in the last section.
 
 To our opinion the approach to approximate optimal analysis, using the direct monitoring of error propagation by the algebraic and analytic properties of external numbers, outgrows the possibilities of Error Analysis, ordinary set-valued analysis and classical asymptotical methods. 
 
 The article has the following structure. In Sections \ref{secion 1c6} and \ref{optimization} we consider some basic notions and properties related to external numbers, flexible functions and near-optimization. The various forms of imprecise limits, approximate continuity and differentiation, together with their behavior under operations, are studied in Section \ref{limits} and Section \ref{section7c6}. Section \ref{section2c8} contains a Fermat Lemma for near-optimization. The remaining part of the article is devoted to flexible functions of several variables, where for reasons of notational and conceptual simplicity we restrict ourselves to functions of two variables. Section \ref{section9c6} starts with a definition of partial derivatives for flexible functions of more variables. In this setting we extend the Fermat Lemma, introduce a total differential and prove a chain rule. Section \ref{section11c6} contains an Implicit Function Theorem under imprecise conditions. All this material is joined in Section \ref{seclag}, to state and prove the Main Theorem on the existence of Lagrange multipliers for the mixed optimization problem of a flexible objective function and a precise constraint. 
 
 The present article is written using the axiomatic approach to Nonstandard Analysis $ IST $ of Nelson \cite{Nelson}, in combination with the theory of external sets of \cite{KanoveiReeken}. Introductions to $ IST $ are for example contained in \cite{Dienerreeb} and \cite{Lyantsekudryk}, and \cite{DinisBerg} contains an introduction to external numbers.

\section{Preliminaries}\label{secion 1c6}
	\subsection{Neutrices and external numbers}\label{neutrices} 
In this subsection we recall the notions of neutrices and external numbers, as well as operations on them. We also recall some background on Nonstandard Analysis and properties which are needed in the sequel.

The axiomatic system $ IST $ postulates the existence of nonstandard elements within infinite standard sets. For instance, the standard set $ \R $ contains infinitesimals and unlimited numbers, larger than any standard real number. Sets of classical set theory $ ZFC $ are called \emph{internal}. Sets of internal elements are often outside $ ZFC $, and then are called \emph{external}. As a consequence of Nelson's Reduction Algorithm \cite{Nelson}, when definable in $ IST $, every external subset of a standard set $ X $ can be expressed in the form $\bigcup\limits_{st(u)\in U} \bigcap\limits_{st(v)\in V} I_{uv}$ where $U, V$ are standard sets and $I : \ U\times V \rightarrow \P(X)$ is an internal set-valued mapping. External sets which reduce to $\bigcup\limits_{st(u)\in U}  I_{u}$ are called \emph{galaxies}, and are called \emph{halos} when they reduce to  $\bigcap\limits_{st(v)\in v}  I_{v}$.
\begin{definition}A \emph{(scalar) neutrix} is an additive convex subgroup of $\R$. An {\em external number}  is the Minkowski-sum of a real number and a neutrix. 
\end{definition}
Classically, the only neutrices are $\{0\}$ and $\R$, but allowing for external sets there are many more neutrices. Examples are $\oslash$, the set of infinitesimals and $\pounds$, the set of limited numbers. Let $\varepsilon\in \R$ be a positive infinitesimal. Other examples of neutrices are $\varepsilon\pounds$, $\varepsilon \oslash $,  $M_{\varepsilon}\equiv\displaystyle\bigcap_{st(n)\in \N}[-\varepsilon^n, \varepsilon^n]=\pounds\varepsilon^{\not \hskip -0.05cm\infty}$ and $\mu_{\varepsilon}\equiv\displaystyle\bigcup_{st(n)\in \N}[-e^{-1/(n\eps)}, e^{-1/(n\eps)}]=\pounds e^{-@/\eps}$, where $ \not \hskip -0.12cm\infty $ is the set of positive unlimited numbers and $ @ $ is the set of positive \emph{appreciable} numbers, i.e. limited numbers which are not infinitesimal. One shows \cite{VDBNAA} that as groups they are not isomorphic, and that every external neutrix is either a halo or a galaxy. A neutrix $ I $ is called \emph{idempotent} if $ II=I $, which is the case of all the examples above with the exception of $\varepsilon\pounds$ and $\varepsilon \oslash $. As is the case for these neutrices, for any neutrix $ N $ there exists $ p\in \R $ such that $ N=pI$, with $ I$ idempotent \cite{Koudjeti Van den Berg}.

External neutrices $ N $ may be seen model for imprecisions of the Sorites type, for they are stable under some shifts, additions and multiplications. In particular they are invariant by multiplication by appreciable numbers, so $@N=N  $. An \emph{absorber} of $N $ is a real number $a $
such that $aN\subset N$ and an \emph{exploder} is a real number $b $
such that $bN\supset N$; the set of absorbers of $ N $ is denoted by $ \oslash_{N} $. In the examples above the infinitesimal number $%
\varepsilon $ is an absorber of $\pounds$ and $\oslash $, and $%
1/\varepsilon $ an exploder of $\pounds$ and $\oslash $. These numbers leave $M_{\varepsilon} $ and $ \mu_{\varepsilon} $ invariant, so in a sense the latter neutrices are more imprecise. 
Neutrices are ordered by inclusion, and if the neutrix $ A$ is contained in the neutrix $ B$, we have $ B=\max\{A,B\} $.

Each external number  has the form $\alpha=a+A=\{a+x|x\in A\}$, where $A$ is called the {\em neutrix part} of $\alpha$, denoted by ${\rm{N}}(\alpha)$, and $a\in \R$ is called a {\em representative}  of $\alpha$. We call $\alpha$ {\em zeroless} if $0\not\in \alpha$,  and {\em neutricial} if $ \alpha={\rm{N}}(\alpha) $.

The collection of all neutrices is not an external set in the sense of \cite{KanoveiReeken}, but a definable class, denoted by $\mathcal{N}$. Also the external  numbers form a class, denoted by $ \E $. 

The rules for addition, subtraction, multiplication and division of external numbers respect the rules of informal Error Analysis. In Definition \ref{defnam} below they are defined formally as  Minkowski operations on sets of real numbers.

\begin{definition}\label{defnam}\rm Let $ a,b\in \R $, $A, B$ be neutrices and $\alpha=a+A, \beta=b+B$ be external numbers.
	\begin{enumerate}[(a)]
		\item \label{cong} $\alpha\pm\beta=a\pm b+A+B= a\pm b+\max\{A,B\}$.		
		\item \label{nhan}  $\alpha \beta=ab+Ab+Ba+AB=ab+\max\{aB, bA, AB\}.$
		\item \label{1/alpha} If $\alpha$ is zeroless, $\frac{1}{\alpha}=\frac{1}{a}+\frac{A}{a^2}.$		
	\end{enumerate}
\end{definition}

If $\alpha$ or $ \beta$ is zeroless, in Definition \ref{defnam}.\ref{nhan} we may neglect the neutrix product $ AB $. An order relation is given as follows.

\begin{definition} Let $ \alpha,\beta \in \E$. We define 
	\begin{equation*}
		\alpha\leq\beta\Leftrightarrow\forall a\in \alpha \exists b\in \beta(a\leq b).
	\end{equation*}
	If $ \alpha \cap \beta =\emptyset $ and $ \alpha\leq\beta $, then $ \forall a\in \alpha \forall b\in \beta(a<b) $ and we write $ \alpha<\beta $.
\end{definition}

In \cite{Koudjeti Van den Berg, dinisberg 2016} it is shown that, with some small adaptations, the relation $ \leq  $ is a total order relation compatible with the operations, while "Generalized Dedekind Completeness" holds for definable subsets of $ \R $. The inverse order relation is given by
\begin{equation*}
	\alpha\geq\beta\Leftrightarrow\forall a\in \alpha \exists b\in \beta(a\geq b),
\end{equation*}
and $ \alpha>\beta $ if $ \forall a\in \alpha \forall b\in \beta(a>b) $.
Clearly $ \alpha<\beta $ implies $ \beta<\alpha $. In general, for external numbers $ \alpha,\beta $ which are not disjointed, the more precise number satisfies both order relations with respect to the less precise number, but not the inverse relation, i.e., if $ \alpha\subseteq\beta $, then $ \alpha\leq \beta$ and $ \alpha\geq \beta$, but $ \beta\leq\alpha $ and $ \beta\geq\alpha $ only if  $ \alpha= \beta$. So both $ \oslash\leq \pounds $ and $ \oslash\geq \pounds $, while $ \pounds\nleq\oslash $ and $\pounds\ngeq\oslash  $. External numbers $ \alpha $ such that $ 0\leq \alpha $ are called \emph{non-negative}. The \emph{absolute value} of an external number $\alpha=a+A$ is defined by $|\alpha|=|a|+A$. Notice that this definition does not depend on the choice of the representative of $\alpha.$

By the close relation to the real numbers, practical calculations with external numbers tend to be quite straightforward. A full list of axioms for the operations on the external numbers has been given in  \cite{dinisberg 2016} and \cite{DinisBerg}, leading to a structure called a \emph{Completely Arithmetical Solid $(CAS)$}. 

Some care is needed with distributivity. Theorem \ref{tcopti2} states that it holds
up to a neutrix. 

\begin{theorem}{\rm \cite{Dinis}}{(Distributivity with correction term)}\label{tcopti2} Let $\alpha, \beta, \gamma=c+C$ be external numbers. Then 
	\begin{equation}\label{dis}
		\alpha \gamma+\beta \gamma=(\alpha+\beta)\gamma+C \alpha+C \beta.
	\end{equation}  
\end{theorem}

Theorem \ref{tcopti3} says that the common distributive law holds under fairly general conditions, i.e. the correction terms figuring in \eqref{dis} may be neglected. To this end we recall the notions of \emph{relative uncertainty} and \emph{oppositeness}.

\begin{definition}{\rm \cite{Koudjeti Van den Berg, Dinis}}\label{dnbruno}
	Let $\alpha=a+A$ and $  \beta=b+B $ be external numbers and $ C $ be a neutrix.
	\begin{enumerate}[(a)]
		\item \label{ru}The \emph{relative uncertainty}  $ {\rm{R}}(\alpha) $  of $\alpha$ is defined by $ {\rm{R}}(\alpha)=A/\alpha $ if $\alpha $ is zeroless, otherwise $ {\rm{R}}(\alpha)=\R $.
		\item\label{opp} $\alpha$ and $  \beta$ are \emph{opposite} with respect to $ C $ if $ (\alpha+\beta)C\subset \max (\alpha C, \beta C). $
	\end{enumerate}
\end{definition}

\begin{theorem}\label{tcopti3}Let $\alpha, \beta, \gamma=c+C$ be external numbers. Then $\alpha \gamma+\beta \gamma=(\alpha+\beta)\gamma$ if and only if    
	${\rm{R}}(\gamma)\subseteq \max( {\rm{R}}(\alpha), {\rm{R}}(\beta))$, or  $\alpha $ and  $\beta $ are not opposite with respect to $ C $. 
\end{theorem}
We see in particular that distributivity may not hold when multiplying with two almost opposite numbers, which is relevant for differentiation. However the subsdistributivity of Theorem \ref{tcopti1} always holds, and in many cases still enables an effective error-analysis. 

\begin{theorem}{(Subdistributivity)}\label{tcopti1}  
	Let $\alpha, \beta, \gamma$ be external numbers. Then $(\alpha+\beta)\gamma\subseteq\alpha \gamma+\beta \gamma$.
\end{theorem}

We end with a proposition listing some useful properties of external numbers. 
\begin{proposition}{\rm {\cite{Koudjeti Van den Berg}}}\label{tcopt} 
	Let $\alpha=a+A$  be a zeroless external number, and $ \gamma $ be an external number, $ B$ be a neutrix and  $n\in \mathbb{N}$ be standard. Then 
	\begin{enumerate}[(a)]
		
		\item \label{tcoptiv}$\alpha B=aB$ and $	\frac{B}{\alpha}=\frac{B}{a}$.
		\item \label{tcoptid}${\rm{N}}(1/\alpha)={\rm{N}}(\alpha)/\alpha^{2} $.
		\item \label{tcoptix}  ${\rm{R}}(\alpha),{\rm{R}}(1/\alpha)\subseteq \oslash$. 	
		\item \label{giao} $\alpha		\cap \oslash \alpha=\emptyset$.	
		\item \label{tcoptp} $ {\rm{N}}(\alpha\gamma )=\alpha {\rm{N}}(\gamma )+ {\rm{N}}(\alpha)\gamma $.	
		\item \label{lei3} If $\alpha $ is limited and is not an absorber of $B$, then $\alpha B=\frac{B}{\alpha}=B.$
	\end{enumerate}
\end{proposition}

\subsection{Some topological notions}\label{section2c6} 

In classical mathematics a neighborhood of a point $a\subseteq\R^{n} $ contains some open ball $ B(a,r) $ with radius $r>0$, and $ a $ is an accumulation point of a set $X\subseteq \R^{n}$ if $ B(a,r)\cap X\setminus \{a\}\neq \emptyset$ for every $r>0$; one may take also products of open intervals. We adapt these notions to external points.

\begin{convention}\label{convXn}
From now on we consider only spaces $ \R^n $ for $ n\in \N $ standard. A subset $X\subseteq \R^n$ is always supposed to be non-empty. 
\end{convention}

We extend the notion of neutrix to convex subgroups $ M\subseteq\R^{n} $; these can be written in the form $ M=M_{1}u_{1}\oplus\cdots \oplus M_{n}u_{n}$, where $M_{1},\dots, M_{n}$ are scalar neutrices and $u_{1},\dots, u_{n}\in\R^{n}$ are orthonormal vectors \cite{Immedecompositionofneutrices}. So an \emph{external point} $ \xi=p+M\in \R^{n} $ is of the form $ (p_{1},\dots, p_{n})+ M_{1}u_{1}\oplus\cdots \oplus M_{n} u_{n}$, with $ p_{i}\in \R $ and $ M_{i} \subseteq\R$ a neutrix for $ 1\leq i\leq n $. Since we are studying limits in more variables, we will always assume that $ (u_{1},\dots, u_{n})=(e_{1},\dots, e_{n}) $ is the canonical basis of $ \R^{n} $.

\begin{definition}\label{defneigh} Let $ n\in \N $, $\xi=p+M=(p_{1},\dots, p_{n})+(M_{1}e_{1}\oplus\cdots\oplus M_{n}e_{n})\in E^{n}$ be an external point and $U\subseteq \R^{n}$ be internal. The set $U$ is said to be an 
	\emph{$M$-neighborhood} of $\xi$ if and only if there exist $ r_{1},\dots, r_{n}\in \R,r_{1}>M_{1},\dots, r_{n}>M_{n} $ such that $ B(p_{1},r_{1})\times\cdots\times B(p_{n},r_{n}) \subseteq U$. Then $\xi$ is called an \emph{$M$-interior point} of $X$ and $ U\setminus\xi $ an \emph{outer $M$-neighborhood} of $\xi$.
	
\end{definition}

If the components $ M_{i}, 1\leq i\leq n $ of $ M $ are different,  neighborhoods in terms of products of intervals can be smaller than open balls around external points, for we have to take a radius larger than the biggest of the $ M_{i} $. It is not true in general that for every internal open set $ U\supset M $ there exists a product of intervals $ B(p_{1},r_{1})\times\cdots\times B(p_{n},r_{n}) \subseteq U$; take for example $ \xi=\pounds\times \oslash $ and $ U=\{(x,y)|-1/x<y<1/x,x\in\R\} $. The property holds if $ M_1,\dots,M_n $ are galaxies \cite{VDBNAA}.

\begin{definition} 
	Let $ n\in \N $, $X\subseteq \R^{n}$ and $p\equiv(p_{1},\dots,p_{n})\in \R^n$. Let $M=M_{1}e_{1}\oplus\cdots\oplus M_{n}e_{n}\subseteq E^{n}$ be a neutrix. We say that $p$ is an \emph{$M$-accumulation point} of $X$ if for all $\delta\equiv(\delta_{1},\dots,\delta_{n}),\delta_{1}>M_{1}, \dots, \delta_{n}>M_{n}$ one has $B(p_{1},\delta_{1})\times\cdots \times B(p_{n},\delta_{n})\cap \left(X\setminus \{p+M\}\right) \not= \emptyset$. 
\end{definition}

\subsection{Flexible functions}

\begin{definition} \rm 
	Let $ n\in \N $ and $X\subseteq \R^n$. A {\emph{flexible function}} is a mapping $F : \   X\rightarrow \E$, such that its graph $ \Gamma(F) \subseteq X\times \R$ is definable in $ IST $.  A function $f:X\rightarrow\R$ such that $f(x)\in F(x)$ for all $x\in X$ is called a \emph{representative} of $F$. The mapping ${\rm{N}}_F :  \ X \longrightarrow \E$ defined by ${\rm{N}}_F(x)={\rm{N}}(F(x))$ for $x\in X$ is called the \emph{neutrix part} of $F$; observe that the neutrix-part is also a flexible function.
\end{definition}

\begin{convention}\label{conv}
	  All scalar neutrices $ M $ will be supposed to be strictly contained in $ \R $. A neutrix in $ \R^n $, where $ n\in \N $, has a decomposition in terms of scalar neutrices, which thus are supposed to be bounded in $ \R $. Let $X\subseteq \R^n$ and $F : \   X\rightarrow \E$ be a flexible function. We will always assume that there exists a (not necessarily internal) representative $f:X\rightarrow\R$ of $ F $, i.e. for all $x\in X$ we have  $F(x)=f(x)+{\rm{N}}_F(x)$. Moreover, by the above for every $ x\in X $ it holds that $ {\rm{N}}_F(x) \neq \R$; this means that $|F(x)|<c  $ for some $ c\in \R $.
\end{convention}

\begin{examples}\label{exflex}
Let $ k\in \N $ be standard, $\alpha_{0},\alpha_{1},\dots, \alpha_{k}  $ be external numbers and $  \alpha_{i}=a_{i}+A_{i} $ for $ 1\leq i\leq k $. A polynomial $ P $ given by $ P(x)=\alpha_{0}+\alpha_{1}x+\cdots +\alpha_{k}x^{k} $ is a flexible function; the ordinary polynomial $ p $ given by $ p(x)=a_{0}+a_{1}x+\cdots +a_{k}x^{k} $ is a representative, and $ P $ can be seen as a perturbation of $ p $, where every coefficient has its individual imprecision $ A_{i} $. Flexible functions $F$ may also be defined by cases, corresponding to external intervals, like 	
\begin{equation}\label{exF}
F(x)=\begin{cases} e^{x}  & \mbox{ if } x\in \pounds\\
	\ln|x| & \mbox{ if } x\not\in \pounds.
\end{cases} 
\end{equation}	
\end{examples}

\section{Optimization problems with a flexible objective function} \label{optimization}

We study optimization problems of the form
\begin{equation}\label{prob}
	\min\limits_{x\in X} F(x), \mbox{  } \max_{x\in X} F(x)
\end{equation}
where $F$ is a  flexible objective function defined on some set $X\subseteq \R^n$. (Nearly) optimal solutions are defined as follows.

\begin{definition} Let $X\subseteq \R^n$ and $F :  \ X \longrightarrow \E$ be a flexible function.  Let $N$ be a neutrix and $a \in X$. 
	\begin{enumerate}[(a)]
		\item The point $a $ is called an \emph{$N$-minimizer}  of the minimization problem $\min\limits_{x\in X}F(x)$ if $F(x)\geq F(a )+N$ for all $x\in X$.
		Then $F(a )$ is called an \emph{$N$-minimum}.
		\item The point $a $ is called an {\it $N$-maximizer} of the maximization problem $\max\limits_{x\in X}F(x)$ if $F(x)\leq F(a )+N$ for all $ x\in X.$
		Then $F(a )$ is called an \emph{$N$-maximum}.
	\end{enumerate}
An $N$-minimal or $N$-maximal value is also called the {\it $N$-optimal value}, and an $N$-minimizer or $N$-maximizer an \emph{$N$-extreme point}. 	
	When $N\subseteq {\rm{N}}(F(a ))$, we may simply call an $N$-maximum a \emph{maximum} and an $N$-minimum a \emph{minimum}. 
\end{definition}

\begin{examples} Let $ \varepsilon\simeq 0, \varepsilon> 0 $. 
	Let $G : \ \R \longrightarrow \E$ be given by  $G(x)=x^2+\varepsilon \pounds $. Then $ G $ admits $ \varepsilon\pounds $-minimizers, and the set of $   \varepsilon\pounds $-minimizers is equal to $ \sqrt{\varepsilon}\pounds $, illustrating that the size of the set near-minimizers tends to be larger than the size of the set of near-minima. Let $H :  \ \R \longrightarrow \E$ be given by  $H(x)=x^2+\oslash x, \ x\in \R$. Though $ H(0)=0 $ is precise, it is not a minimizer for the optimization problem $\min\limits_{x\in \R} H(x)$, because for all $ x\simeq 0 $ it holds that $  x^{2}+\oslash x = \oslash x $, and the neutrix $ \oslash x $ contains negative numbers. However, $ 0 $ is an $\oslash$-minimizer, because $   \oslash x \geq \oslash$ for $ x\simeq 0 $ and $  x^{2}+\oslash x > x^2/2 >\oslash$ for $ x\notin \oslash $. In fact,  every $a \in \oslash$ is an $\oslash$-minimizer of $F$ on $\R$. 
	
	In general, it may be expected that near-extreme values of the polynomials $ P $ mentioned in Examples \ref{exflex} come as perturbations of extreme values of the polynomials $ p$. The optimization problem for the flexible function $ F $ of \eqref{exF} is more intricate, even if $ F $ is differentiable at every $ x\in \R $, decreasing for $ x<\pounds $ and increasing for $ x\geq \pounds $; when the variable decreases in $ \pounds $, the values of $ e^{x} $ approach the infinitesimals, but the weak infimum $ \oslash $ in the sense of Generalized Dedekind Completeness is never attained.
\end{examples}

The above examples suggest that the set of near-optimal solutions is often external, of the form of an external point $ a+M $. To define local near-optimal solutions we must consider behavior on the neighborhoods of such external points, as given by Definition \ref{defneigh}.

\begin{definition}\label{M,N local optimal solution}\rm  Let $X\subseteq \R^n$, $F :  \ X \rightarrow \E$ be a flexible function and $a\in X$. Let $M\subseteq \R^{n}$ and $N\subseteq \R$ be neutrices.  
	\begin{enumerate}[(a)]
		\item The point $a$ is called an {\it $M$-local $N$-minimizer} of the problem $\min\limits_{x\in X} F(x)$ if there exists an $ M $-neighborhood $ U $ of $ a+M $ such that $ a $ is an  $N$-minimizer of the problem $\min\limits_{x\in U} F(x)$.
		\item The point $a$ is called an {\it $M$-local $N$-maximizer} of the problem $\max\limits_{x\in X} F(x)$ if there exists an $ M $-neighborhood $ U $ of $ a+M $ such that $ a $ is an  $N$-maximizer of the problem $\max\limits_{x\in U} F(x)$.
	\end{enumerate}
	In particular, if $N={\rm{N}}_F(a)$, for  minimization problems  we call $a$ an {\it $M$-local minimizer} and for  maximization problems  we call $a$ an {\it $M$-local maximizer}. 
	
\end{definition}

\begin{example}\rm  Let $F :  \ \R \rightarrow \E$ be the flexible function defined by $F(x)=x^3-3x+1+\oslash x$. It is easy to verify that every $x\simeq  -1$ is an $\oslash$-local maximizer  and every $x\simeq 1$ is an $\oslash$-local minimizer of $F$, but that these points are neither $\pounds$-local minimizers of $F$ nor $\pounds$-local maximizers of $F$.
\end{example}

Propositions \ref{nghiem neutrix nho} shows that an $M,N$-local optimal solution is also an $M',N'$-local optimal solution for $M'\subseteq M$ and $N'\supseteq N$. As a consequence, in practice we tend to determine $M,N$-local optimal solutions with the largest possible $M$ and the smallest possible $ N $, i.e. the $ M $-minimizers and $ M $-maximizers of Definition \ref{M,N local optimal solution}.  

\begin{proposition}\label{nghiem neutrix nho} Let $F :  \ X \subseteq \R^n \longrightarrow \E$ be a flexible function and $ M,M',N,N' $ be neutrices such that $ M'\subseteq M\subset\R^{n}$ and $N\subseteq N'\subset\R$. 
	\begin{enumerate}[(a)]
	\item \label{nghiem neutrix nho1} Assume that  $a\in X$ is an $M$-local $N$-optimal solution of  \eqref{prob}. Then $a$ is an $M'$-local $N$-optimal solution of this problem.
	\item \label{nghiem neutrix lon}
	Assume that  $a\in X$ is an $M$-local $N$-optimal solution of  \eqref{prob}. Then $a$ is an $M$-local $N'$-optimal solution of this problem.
\end{enumerate}
\end{proposition}
\begin{proof}Without loss of generality we assume that $a$ is an $M$-local $N$-minimizer of \eqref{prob}.
	
	\eqref{nghiem neutrix nho1} 
	There exists an $ M $-neighborhood $ U\subseteq X $ of $a+M$ such that $F(x)\geq F(a)+N$ for all $x\in U$. Because  of $a+M\supseteq a+M'$ the point $a$ is also an $ M',N $-minimizer of $ F $ on $ U $, i.e. an $M'$-local $N$-minimal solution of the problem \eqref{prob}.

	\eqref{nghiem neutrix lon} Because $a$ is an $N,M$-local minimizer of $F$ on $X$, there exists an $ M $-neighborhood $ U\subseteq X $ of $a+M$ such that $F(x)\geq F(a)+N$ for all $x\in U$. Now $ N\geq N' $, so by transitivity we have $F(x)\geq F(a)+N'$ for all $x\in U$. Hence $ a $ is an $M$-local $N'$-minimal solution of the problem \eqref{prob}.
\end{proof}

\section{Limits of flexible functions}\label{limits}
In \cite{Sequences} limits of sequences up to a neutrix $ N $ were introduced. Some notions will be adapted to study the behaviour of a flexible function when the variables  approach an external point, in particular for minimal $ N $. The variables may stay outside the external point (outer limits) or enter it (inner limits), and we consider also a mixed form.
\begin{definition} \label{dnouterlimits}
	Let $ n\in\N $, $X\subseteq \R^n$ and $F: X\rightarrow \E$ be a flexible function. Let $M=(M_1,M_2,\dots, M_n)$ be a neutrix vector, $N$ be neutrix  and  $a=(a_1, a_2,\dots, a_n)$ be an $M$-accumulation point of $X$. An external number $\alpha$ is called an \emph{$M\times N$-outer limit} of $F$ at $a$,  if for all $\eps>N$ and for every $ i $ with $1\leq i\leq n$ there exists $\delta_i>M_i$  such that for all $x=(x_1, x_2,\dots, x_n)\in X$ 
	$$ \forall i, 1\leq i\leq n,  M_{i}< |x_i-a_i|<\delta_i\Rightarrow |F(x_{1},x_2,\dots, x_{n})-\alpha|<\eps .$$
	Then we write 
	$$N\mbox{-}\lim\limits_{\substack{x\to a+M}} F(x)=\alpha.$$
	The \emph{$M$-outer limit} of $F$ at $a$ is defined by 
	\begin{equation}\label{olim}
	\lim\limits_{\substack{x\to a+M}} F(x)=\bigcap\limits_{\alpha={\rm{N}}(\alpha)\mbox{-}\lim\limits_{\substack{x\rightarrow a+M}} F(x)}\alpha. 
	\end{equation} 	
	We call ${\rm{N}}_{F}^{(M)}(a)\equiv {\rm{N}}\Big( \lim\limits_{x\to a+M} F(x)\Big) $ the \emph{$ M $-outer limit neutrix} of $ F $ at $ a $. 
\end{definition}
If $ F $ is bounded, the $ N $-outer limit exist for $ N $ larger than the lowest upper bound for $ |F| $ (in the sense of Generalized Dedekind Completeness). The limit is not unique, for if $\alpha$ is an $M\times N$-outer limit of $F$ at $a$, every external number $\alpha'  $ with   $\alpha'\subseteq\alpha+N$ is an $M\times N$-outer limit of $F$ at $a$. Note that also the neutrix $ N $ is not unique since if $N\mbox{-}\lim\limits_{x\to a+M} F(x)$ exists, also $N'\mbox{-}\lim\limits_{x\to a+M} F(x)$ is well-defined for every neutrix $ N'\supseteq N $. We get uniqueness if we let the neutrix $ N $ be the weak infimum of the neutrices $ N' $, and $ {\rm{N}}(\alpha)=N $, obtaining maximal precision for the limit. It follows from
\cite[Th. 5.1.2]{DinisBerg} that the weak infimum corresponds to a minimum. Hence ${\rm{N}}^{(M)}(a)  $ is the minimal neutrix for which the outer limit of $ F $ at $a$ exists, i.e., if $ a\in \bigcap\limits_{\alpha={\rm{N}}(\alpha)-\lim\limits_{x\to a+M} F(x)} \alpha $, the outer limit is an $ M\times {\rm{N}}_{F}^{(M)}(a)$-limit and we have
\begin{equation}
	\lim\limits_{x\to a+M} F(x)={\rm{N}}_{F}^{(M)}(a)\mbox{-}\lim\limits_{x\to a+M} F(x)=a+{\rm{N}}_{F}^{(M)}(a).
\end{equation} 
\begin{definition}\label{FMOcont}
	Let $X\subseteq \R^{n}$, $F :  \ X\longrightarrow \E$ be a flexible function and $ \alpha$ be an external number. Let $ M=(M_{1},M_2, \dots, M_{n})$ be a neutrix vector and $N$ be a neutrix. Assume that $a=(a_{1},a_2,\dots,a_{n})$ is an $M$-accumulation point of $X$. The function $F$ is said to be \emph{$M\times N$-outer continuous} at $a$ if $F$ is defined at $a$ and  $N\mbox{-}\lim\limits_{x\to a+M}=F(a)+N$. If $\lim\limits_{x\to a+M}F(x)=F(a)$ we say that $F$ is \emph{$M$-outer continuous} at $a$. 
\end{definition}
Observe that in the case of $ M $-outer continuity we have ${\rm{N}}_{F}^{(M)}(a)={\rm{N}}_{F}(a)$.

For inner limits the limit behaviour of the function outside the external point $ a+M $ must persist when going inside. 

\begin{definition}\label{FMIlimit} Let $ n\in\N $ be standard. Let $X\subseteq \R^n$ and $F: X\rightarrow \E$ be a flexible function. Let $M=(M_1,M_2,\dots, M_n)$ be a neutrix vector, $N$ be neutrix  and  $a=(a_1,a_2, \dots, a_n)$ be an $M$-accumulation point of $X$. An external number $\alpha$ is called an \emph{$M\times N$-inner limit} of $F$ at $a$  if for all $\eps>N$ and for every $ i $ with $1\leq i\leq n$ there exists $\delta_i>M_i$  such that for all $x=(x_1,x_2,\dots, x_n)\in X$ 
$$ \forall i, 1\leq i\leq n,  0< |x_i-a_i|<\delta_i\Rightarrow |F(x_{1},x_2,\dots, x_{n})-\alpha|<\eps .$$
 Then we write 
	$$N\mbox{-}\lim\limits_{\substack{x\twoheadrightarrow a+M}} F(x)=\alpha.$$
The \emph{$M$-inner limit} of $F$ at $a$ is defined by 
	$$\lim\limits_{\substack{x\twoheadrightarrow a+M}} F(x)=\bigcap\limits_{\alpha={\rm{N}}(\alpha)\mbox{-}\lim\limits_{\substack{x\twoheadrightarrow a+M}} F(x)}\alpha.$$ 
	The flexible function $ F $ is called \emph{$ M $-inner continuous} at $ a $ if $ F $ is defined at $ a $ and $ F(a)= \lim\limits_{x\twoheadrightarrow a+M} F(x) $.
	
\end{definition}

Theorem \ref{inclusion of both limits} states that outer limits are included in inner limits, and the imprecision of $ M $-inner limits grows with $ M $. 

\begin{theorem}\label{inclusion of both limits} 
	Let $X\subseteq \R^{n}$ and $F: X\rightarrow \E$ be a flexible function. Let $M=(M_{1},M_2,\dots, M_{n})\subset\R^{n}$ be a neutrix and  $a$ be an $M$-accumulation point of $X$. Assume $ \lim\limits_{\substack{x\twoheadrightarrow a+M}} F(x) $ is well-defined.
	\begin{enumerate}[(a)]
	\item \label{ininout} $\lim\limits_{x\to a+M} F(x)\subseteq\lim\limits_{\substack{x\twoheadrightarrow a+M}} F(x)$. 
	\item \label{inin} If $ M'=(M'_{1},M_2',\dots, M'_{n})\subseteq M $ is a neutrix, then  $\lim\limits_{x\twoheadrightarrow a+M'} F(x)\subseteq\lim\limits_{\substack{x\twoheadrightarrow a+M}} F(x)$. 
\item \label{limin} If $ F\equiv f: X\rightarrow \R^{n}$ is internal and $ \lim\limits_{x\to a} f(x) $ is well-defined, then  $\lim\limits_{x\to a} f(x)\in\lim\limits_{\substack{x\twoheadrightarrow a+M}} f(x)$, and $f(a)\in\lim\limits_{\substack{x\twoheadrightarrow a+M}} f(x)$, if in addition $ f $ is continuous at $ a $. 
	\end{enumerate}
\end{theorem}
\begin{proof}Let $\alpha= \lim\limits_{x\twoheadrightarrow a+M} F(x)$,  $\eps>{\rm{N}}(\alpha)$ and $ a=(a_1,a_2, \dots, a_n)$, with $a_1,a_2, \dots, a_n\in \R$. There exist  $\delta_1>M_1,\delta_2>M_2,\dots,\delta_n>M_n$,  such that for all $x=(x_1,x_2, \dots, x_n)\in X$ with $0< |x_i-a_i|<\delta_i$ for $ 1\leq i\leq n$
	\begin{equation}\label{Feps}
	|F(x_1,x_2, \dots, x_n)-\alpha|<\eps
	\end{equation}
	\begin{enumerate}[(a)]
		\item In particular \eqref{Feps} holds for $ M_{i}< |x_i-a_i|<\delta_i,1\leq i\leq n $. So $ {\rm{N}}(\alpha)\mbox{-}\lim\limits_{x\to a+M} F(x)= \alpha$, hence  $\lim\limits_{x\to a+M} F(x)\subseteq \alpha$ by minimality of ${\rm{N}}\left(\lim\limits_{x\to a+M} F(x)\right)$.
		\item If $ M'_{i}<\delta'_{i} \leq \delta_i$ for $ 1\leq i\leq n$, formula \eqref{Feps} holds for $ 0< |x_i-a_i|<\delta'_i $. So $ {\rm{N}}(\alpha)\mbox{-}\lim\limits_{x\to a+M'} F(x)= \alpha$, hence  $\lim\limits_{x\to a+M'} F(x)\subseteq \alpha$ by minimality of ${\rm{N}}\left(\lim\limits_{x\to a+M'} F(x)\right)$.
		\item The properties follow directly from Part \eqref{inin}, taking $ M'=\{0\} $.
	\end{enumerate}

\end{proof}

Generally speaking, limits of flexible functions respect the algebraic operations with inclusion, with the exception of the product of limits of functions tending to a neutrix.  Theorem \ref{operationsNconv} gives the rules for outer limits for functions of one variable.
\begin{theorem} \label{operationsNconv}
	Let $ X\subseteq\R $, $M$ be a neutrix and $a\in \R$. Assume that $ a+M $ is an $ M$-accumulation point of $ X $. Let $ F,G:X\rightarrow \R $ be flexible functions. Let $\alpha ,\beta \in \mathbb{E}$ be such that $ \lim\limits_{x\to a+M} F(x)=\alpha $ and $ \lim\limits_{x\to a+M} F(x)=\beta $. Then
	
	\begin{enumerate}[(a)]
		\item \label{addition} $ {\rm{N}}\left(  \lim\limits_{x\to a+M} (F+G)(x)\right) \subseteq {\rm{N}}\left( \lim\limits_{x\to a+M} F(x)\right) +{\rm{N}}\left(  \lim\limits_{x\to a+M} G(x)\right)  $ and $ \lim\limits_{x\to a+M} F(x)+  G(x)\subseteq	 \alpha +\beta $.						
		\\
		\item \label{Theorem product}
		
		\begin{enumerate}[(i)]
		\item \label{zorz}If $ \alpha $ or $ \beta $ is zeroless,
		
	 $ {\rm{N}}\left(  \lim\limits_{x\to a+M} (FG)(x)\right) \subseteq\alpha {\rm{N}}\left( \lim\limits_{x\to a+M} G(x)\right) +\beta {\rm{N}}\left(  \lim\limits_{x\to a+M} F(x)\right) $ and $ \lim\limits_{x\to a+M} (FG)(x)\subseteq	 \alpha \beta $.
		\item \label{nan} If $ \alpha \equiv A$ and $ \beta\equiv B $ are neutricial, 
		$ \lim\limits_{x\to a+M} (FG)(x)\subseteq {\rm{N}}\left(  \lim\limits_{x\to a+M} (FG)(x)\right) \subseteq qA+pB$, where $ p,q\in \R $ are such that $ A=pI, B=qJ $, with $ I,J $ idempotent neutrices.
			\end{enumerate}		
		\item \label{prodprecise} Assume (i) $ c\in \R $, (ii) $ c\in \E $ is zeroless or (iii)  $ \alpha $ is zeroless. Then $ {\rm{N}}\left(  \lim\limits_{x\to a+M} (cF)(x)\right) =c{\rm{N}}\left( \lim\limits_{x\to a+M} F(x)\right) $ and $ \lim\limits_{x\to a+M} cF(x)=	c \alpha $.
		\item \label{division} If $F$ is zeroless in some outer $ M $-neighborhoud of $ a $, and $\alpha$ is zeroless, then ${\rm{N}}\left( \lim\limits_{x\to a+M} \frac{1}{F(x)}\right) =\frac{{\rm{N}}\left( \lim\limits_{x\to a+M}F(x)\right)}{\alpha^{2}}$ and  $\lim\limits_{x\to a+M} \frac{1}{F(x)}=	\frac{1}{\alpha}$.
	\end{enumerate}
	\end{theorem}

\begin{proof} Let $ \alpha=s+A, \beta=t+B $, with ${\rm{N}}(\alpha)=A$ and ${\rm{N}}(\beta) =B$. Observe that $ F,G $ are bounded in some outer $ M $-neighborhood of $ a $, hence also $ F+G,FG $ and $ 1/F $, which means that their $ N $-limits exist for sufficiently large neutrices $ N $.
	
	\eqref{addition} Let  $C={\rm{N}}\left( \lim\limits_{x\to a+M} (F+G)(x)\right)$. Let $\eps>A+B$. Then $\eps/2>A$ and $\eps/2>B$. There exist $\delta_1,\delta_2>M$ such that $\left| F(x)-\alpha\right|<\eps/2$ for all $x, M<|x-a|<\delta_{1}$ and  $\left| G(x)-\beta\right|<\eps/2$ for all $x, M<|x-a|<\delta_{2}$. Let $\delta=\min\{\delta_1, \delta_2\}$. Then for all $x, M<|x-a|<\delta$ one has 
	\begin{equation}\label{addition1}   \left|(F+G)(x)-\alpha-\beta\right|\leq \left|F(x)-\alpha\right|+\left|G(x)-\beta\right|<\eps/2+\eps/2=\eps.
		\end{equation} 
	Hence $ (A+B)\mbox{-}\lim\limits_{x\to a+M} (F+G)(x)=\alpha +\beta $. By minimality of $C$  we conclude that $C\subseteq A+B$ and $\lim\limits_{x\to a+M} F(x)+  G(x)\subseteq	 \alpha +\beta$.
	
\eqref{Theorem product}(i) 
 Let $D\equiv {\rm{N}}\left(  \lim\limits_{x\to a+M} (FG)(x)\right)$.
	Without loss of generality, we assume that $\alpha$ is zeroless.  Then  
	\begin{equation}\label{prod0}
				\left| (FG)(x)-\alpha\beta \right|  \leq  \left|G(x)\right|\left|F(x)-\alpha\right|+\left|\alpha\right|\left|G(x)-\beta\right|
			\leq  \left|G(x)\right|\left|F(x)-\alpha\right|+\left|2s\right|  \left|G(x)-\beta\right|.
	\end{equation}
Let  $K={\rm{N}}(\alpha\beta)=\alpha B+\beta A=s B+t A$. Let $\eps>K$. Then $\eps>s B$, which implies that $\frac{\eps}{4|s|}>B$. Then there exists $\delta_1>M$ such that for all $x, M<|x-a|<\delta_1$
	\begin{equation}\label{prod1}
		\left|G(x)-\beta\right|<\frac{\eps}{4|s|}.
	\end{equation}	
	 
	 If  $\beta$ is zeroless, then $|t|>B$, so there exists $\delta_2>M$ such that   $|G(x)-\beta|<|t|$ for all $x, M<|x-a|<\delta_2$. Hence
	  \begin{equation}\label{prod2}
		|G(x)|<2|t|
	\end{equation} for all $x, M<|x-a|<\delta_2$. Also, it follows from $\eps>K$  that $\frac{\eps}{4|t|}>A$. There exists $\delta_3>M$ such that 
	\begin{equation}\label{prod3} 
		\left|F(x)-\alpha\right|<\frac{\eps}{4|t|}
	\end{equation} for all $x, M<|x-a|<\delta_3$. Let $\delta\equiv \min\{\delta_1,\delta_2, \delta_3\}$. By \eqref{prod0}-\eqref{prod3} one has 
\begin{align*}
	\left| (FG)(x)-\alpha\beta \right| 
	<2|t|\frac{\eps}{4|t|} +\left|2s\right| \frac{\eps}{4|s|}=\eps.  
\end{align*}
for all $x, M<|x-a|<\delta.$

 If $\beta=B$ is a neutrix 
there exists $\delta_4>M$ such that for all $x, M<|x-a|<\delta_4$ 
\begin{equation}\label{prod4}
	\left|G(x)\right|<\frac{\eps}{2|s|}.
\end{equation}
Also $|s|>A$, so there exists $\delta_5>M$ such that 
\begin{equation}\label{prod5} 
	 \left|F(x)-\alpha\right|< |s| 
	\end{equation}  for all $x, M<|x-a|<\delta_5$. Let $\delta' \equiv  \min\{\delta_1, \delta_4, \delta_5\}$. Combining \eqref{prod0}, \eqref{prod1}, \eqref{prod4} and \eqref{prod5} one obtains
\begin{align*}
	\left| (FG)(x)-\alpha\beta \right| 
	<\frac{\eps}{2|s|}|s| +\left|2s\right| \frac{\eps}{4|s|}=\eps.  
\end{align*}
for all $x, M<\left|x-a\right|<\delta'$.

In both cases we conclude that $K\mbox{-}\lim\limits_{x\to a+M} (FG)(x)=\alpha\beta$. Then it follows from the minimality of $D$ that $D\subseteq \alpha B+\beta A$ and  $\lim\limits_{x\to a+M} (FG)(x)\subseteq \alpha\beta$.

(ii) Let  $\eps>pB+qA\equiv K$. Then $\eps>qA$ and $\eps>pB$. As a result, $\eps>pqI$ and $\eps>pqJ$. So $\eps/pq>I$ and $\eps/pq>J$. Because $I, J$ are idempotent, it follows that $\sqrt{\frac{\eps}{pq}}>I$ and $\sqrt{\frac{\eps}{pq}}>J$. Hence  $\eps_1\equiv\sqrt{\frac{p\eps}{q}}=p\sqrt{\frac{\eps}{pq}}>pI=A$ and   and $\eps_2\equiv\sqrt{\frac{q\eps}{p}}=q\sqrt{\frac{\eps}{pq}}>qJ=B$. There exist $\delta_1, \delta_2>M$ such that $\left|F(x)\right|<\eps_1$ for all $x, M<\left|x-a\right|<\delta_1$ and $\left|G(x)\right|<\eps_2$ for all $x, M<\left|x-a\right|<\delta_2$. Let $\delta=\min\{\delta_1, \delta_2\}$. Then $\left|(FG)(x)\right|=\left|F(x)\right|\left|G(x)\right|<\eps_1\cdot \eps_2=\eps$. So $K\mbox{-}\lim\limits_{x\to a+M} (FG)(x)=K$. By minimality of  $D$ we conclude that $\lim\limits_{x\to a+M} (FG)(x)\subseteq K.$

\eqref{prodprecise}. The proof is obvious in case $ c=0 $. The remaining cases follow from \eqref{Theorem product}.

\eqref{division}. Let $ F(x)=f(x)+K(x)$ for all $x\in X$, with $ f $ real-valued and $ K $ a neutrix-function. Let $\eps>A/s^2$. Then $\eps s^2/8>A$.  So there exists $\delta_1>M $ such that for all $x\in X, M<\left|x-a\right|<\delta_1$ one has $|F(x)-\alpha|<s^2\eps/8$. Now both  $ |f(x)-s|\leq |F(x)-\alpha|$ and $ K(x)\leq |F(x)-\alpha|$, so  $|f(x)-s|<s^2\eps/8$ and  $K(x)<s^2\eps/8$.  Since $|s|>A$, there exists $M<\delta_2$ such that for all $x, M<\left|x-a\right| <\delta_2$ one has $|f(x)-s|\leq |F(x)-\alpha|<|s|/2$, which implies that for all  $x, M<\left|x-a\right| <\delta_2$
\begin{equation}\label{estimate f} 
	|s|/2<|f(x)|<2|s|.
\end{equation}  
Then for all $x, M<\left|x-a\right|<\delta_2 $ one has $s^2/4<f^2(x)<4s^2$.  Let $\delta=\min\{\delta_1, \delta_2\}$. Then for all $x, M<\left|x-a\right|<\delta$ one has 
\begin{alignat*}{2}
	\Big|\frac{1}{F(x)}-\frac{1}{\alpha}\Big| =& \Big|\frac{1}{f(x)}+\frac{K(x)}{f^2(x)}-\frac{1}{s}+\frac{A}{s^2}\Big|
	\leq	\Big|\frac{s-f(x)}{f(x)s}\Big|+\frac{K(x)}{f^2(x)}+\frac{A}{s^2}\\
	\leq  &	\Big|\frac{s-f(x)}{s^2/2}\Big|+\frac{K(x)}{s^2/4}+\frac{A}{s^2}
	< \frac{\eps}{4}+\frac{\eps}{2}+\frac{\eps}{8} <\eps.
\end{alignat*} 
Hence $N\mbox{-}\lim\limits_{x\to a+M} \frac{1}{F(x)}=\frac{1}{\alpha}$ with $N=\frac{A}{s^2}$.

Now we prove that $N$ is minimal. Suppose that $\lim\limits_{x\to a+M}\frac{1}{F(x)}=\frac{1}{s}+C$ with $C \subset N$. We show that $\lim\limits_{x\to a+M} F(x)=s+s^2C$.  Let $\eta>s^2C$. Then $\eta_1\equiv\eta/2s^2>C$.  There exists $ \delta_3 $ with $M<\delta_3\leq \delta_2$ such that 
$$\left|\frac{1}{F(x)}-\frac{1}{s}+C\right|<\eta_1$$ for all $x\in X, M<|x-a|<\delta_3$, and then by \eqref{estimate f} also 
\begin{align*}\label{divisionpro1} 
	\left| 
	\frac{F(x)-s+s^2C}{2s^2}\right|\leq &	\left| 
	\frac{F(x)-s}{s f(x)}+C\right|\leq  \left| \frac{1}{s}-\frac{1}{f(x)}+\frac{K(x)}{f^2(x)/2}+C\right|\\
	=& \left|\frac{1}{f(x)}+\frac{K(x)}{f^2(x)}-\frac{1}{s}+C\right|=\left|  \frac{1}{F(x)}-\frac{1}{s}+C \right|  <\eta_1.
\end{align*} 
This implies that $\left|F(x)-s+s^2C\right| <\eta$ for all $x, M<\left|x-a\right|<\delta_3$.  This means that  $\lim\limits_{x\to a+M}F(x)=s+s^2C$ with $s^2C\subset A$, which is a contradiction to the minimality of $A={\rm{N}}\left(\lim\limits_{x\to a+M}F(x) \right) $.

We conclude that $\lim\limits_{x\to a+M}\frac{1}{x}=\frac{1}{\alpha}$. 
	\end{proof}

We present now a sort of "chain rule", indicating that change of variables for outer limits also has the effect of an inclusion.
\begin{theorem}\label{xylimits}
	Let $ X\subseteq\R $, $M$ be a neutrix and $a\in \R$. Assume that $ a+M $ is an $ M$-accumulation point of $ X $. Let $ \varphi:X\rightarrow \R$, $ b \in \R$ and $ N$ be a neutrix such that $ \lim\limits_{x\to a+M} \varphi(x)=b+N$ and $ \varphi(x)\cap (b+N) =\emptyset $ holds on some  $ M $-outer neighborhood $ U $ of $ a $. Let $G:\varphi(X)\rightarrow \R $ be a flexible function. Assume that $ \lim\limits_{y\to b+N}G(y)=\gamma \in \E$. Then
	$\lim\limits_{x\to a+M} G\circ\varphi(x)\subseteq\gamma$.
\end{theorem}
\begin{proof}Let $\gamma=c+C$ and $\eta>C$. There exists $ \varepsilon>N $ such that whenever $y$ satisfies  $ N<\left|y\right|<\varepsilon$, it holds that $\left| G(y)-\gamma\right|<\eta $. There exists $\delta>M$ such that for all $x, M<\left|x-a\right|<\delta$ one has $\left|\varphi(x)-b\right|\leq\left|\varphi(x)-(b+N)\right|<\eps$, we may assume that   $(a-\delta, a+\delta)\subseteq U$. Then for all $x, M<\left|x-a\right|<\delta$ one has $\left|\varphi(x)-b\right|<\varepsilon$. It follows that for all $x, M<\left|x-a\right|<\delta$ one has $\left|G(\varphi(x))-\gamma\right|<\eta$. So $C\mbox{-}\lim\limits_{x\to a+M} G\circ \varphi(x)=\gamma. $ Hence $\lim\limits_{x\to a+M} G\circ \varphi(x)\subseteq \gamma.$
	
\end{proof}	
In Section \ref{section7c6} we define a sort of a total differential for a flexible function of two variables, where we need to deal with an error function which takes values outside the accumulation point in one variable, but inside for the other variable. To this end we define a mixed inner and outer limit. 
\begin{definition}\label{milim} Let $M=(M_1, M_2)$ be a neutrix vector and $N$ be a (scalar) neutrix. 
	Let $X\subseteq \R^2$ and $F:X\ \to \E$ be a flexible function, $(a, b)$ be both an accumulation point of $X$ and an $ M $-accumulation point of $X$,   and $ \alpha\in \E $. We define the \emph{$N\times M $-mixed limit}  of $ F $ at $(a, b)$ by
	\begin{equation}\label{intout}	
	\begin{array}{ll}
		&	N\mbox{-}\lim\limits_{\substack{x \to a+M_1\\y\twoheadrightarrow b+M_2}} F(x, y)=\alpha \nonumber \\
		\Leftrightarrow  & \forall \eps>N, \exists \delta_1>M_1, \delta_2>M_2\forall (x, y)\in X( M_{1}<|x-a|<\delta_{1}, |y-b|<\delta_2\Rightarrow |F(x, y)-\alpha|<\eps )\\	
		&	N\mbox{-}\lim\limits_{\substack{x \twoheadrightarrow a+M_1\\y\to b+M_2}} F(x, y)=\beta \nonumber \\
		\Leftrightarrow & \forall \eps>N, \exists \delta_1>M_1, \delta_2>M_2\forall (x, y)\in X( |x-a|<\delta_{1}, M_{2}<|y-b|<\delta_2\Rightarrow |F(x, y)-\alpha|<\eps )\nonumber.
	\end{array}
\end{equation}
	In addition we define the \emph{$M  $-mixed limits}
	\begin{equation}\label{outin}	
	\lim\limits_{\substack{x\to a+M_1\\y\twoheadrightarrow b+M_2}} F(x, y)\equiv\bigcap\limits_{\alpha={\rm{N}}(\alpha)\mbox{-}\lim\limits_{\substack{x\to a+M_1\\y\twoheadrightarrow b+M_2}} F(x, y)}\alpha , \hspace{1cm} \lim\limits_{\substack{x \twoheadrightarrow a+M_1\\y \to b+M_2}} F(x, y)\equiv\bigcap\limits_{\beta={\rm{N}}(\beta)\mbox{-}\lim\limits_{\substack{x \twoheadrightarrow a+M_1\\y \to b+M_2}} F(x, y)}\beta.\end{equation}

\end{definition}
	N.B. In Definition \ref{milim} we allow us an abuse of notation, assuming that the functions satisfy the criterium for inner limits also at the limit point.
\begin{remark}\label{oim}
One verifies in a straightforward way that, mutatis mutandis, Theorems \ref{inclusion of both limits}, \ref{operationsNconv} and \ref{xylimits} continue to hold for inner limits, functions of more variables, and mixed limits; in addition in Theorem \ref{xylimits} the type of the limits for the function $ G $ and $ \varphi $ may be different.
\end{remark}
\section{Neutrix-derivatives}\label{section7c6} 

When $h$ approaches $0$, the neutrix part of the expression $\frac{F(a+h)-F(a)}{h}$, in general,  approaches $\R$.   For example, $\frac{\oslash}{h}$ tends to $ \R$ when $h$ approaches $0$. However outer limits with respect to a neutrix $ M\supset\{0\} $ may very well be bounded, and enable to define an $M\times N$-derivative of $ F $, and an $M$-derivative minimizing $ N $. We study the behavior under operations and present a Chain Rule.

\begin{definition}\label{defderNM}Let $M, N$ be neutrices,  $F : \ X\subseteq \R \rightarrow \E$ be a flexible function and $a\in X$ be an $M$-interior point of $X$.  The flexible function $F$ is called \emph{$M\times N$-differentiable}  at $a$ if the $M\times N$-outer limit of the fraction $\frac{F(x)-F(a)}{x-a}$ exists. Then this $M\times N$-outer limit is called the \emph{$M\times N$-derivative}   of $F$ at $a$ and denoted by $\frac{d_N F}{d_M x}(a)$. So
	\begin{equation}\label{defdifNM}
		\frac{d_N F}{d_M x}(a)=N\mbox{-}\lim\limits_{x\to a+M} \frac{F(x)-F(a)}{x-a}.
	\end{equation}
	In case  $N$ is minimal, i.e. if it is the $ M $-limit neutrix of $ \frac{F(x)-F(a)}{x-a} $ at $a$, the limit \eqref{defdifNM} becomes a limit in the sense of \eqref{olim}, and we call 
	\begin{equation}\label{defdifM}
		\frac{d F}{d_M x}(a)\equiv  D_{M}F(a)=\lim\limits_{x\to a+M}\frac{F(x)-F(a)}{x-a}
	\end{equation} 
	the \emph{$M$-derivative} of $ F $ at $ a $. 
\end{definition}

Let $ f : \ X \rightarrow \R $ be a representative of $ F $. In \eqref{defdifM}, put $ x=a+h $. Then 
\begin{equation}\label{decdifM}
	\frac{d F}{d_M x}(a)= D_{M}f(a)+\lim\limits_{h\to M}\frac{{\rm{N}}(F(a))}{h}+\lim\limits_{h\to M}\frac{{\rm{N}}(F(a+h))}{h}.
\end{equation} 
We see that the neutrix-derivative contains two "singular" neutrix terms $\lim\limits_{h\to M}\frac{{\rm{N}}(F(a))}{h}  $ and $ \lim\limits_{h\to M}\frac{{\rm{N}}(F(a+h))}{h} $. We avoid explosion of these terms under the following stability condition. 
\begin{definition}\label{neutder}
A scalar neutrix $ N $ is \emph{stable} for a scalar neutrix $ M $ if $ N $ contains all absorbers of $ M $. 
\end{definition}
Observe that in this case  $ \frac{N}{h} \subseteq N$ for $ h\notin M $, hence
$
\lim\limits_{h\to M}\frac{N}{h}\subseteq\lim\limits_{h\to M}N=N$. As a consequence, the two neutrix terms in \eqref{decdifM} do not exceed ${\rm{N}}(F(a))  $, if $ {\rm{N}}(F(x)) $ is stable for $ M $ in some $ M $-neighborhood of $ a $, and in addition $ F $ is $ M$-outer continuous at $ a $.

Here are some examples. Let $ f:\R\rightarrow\R $ be standard of class $ C^{2} $ in a  neighbourhood of some standard $ a\in \R $, with $ f''(a)\neq 0 $. Then 

\begin{equation}\label{fc2}
D_{\oslash}f(a)= f'(a)+\oslash.
\end{equation}
 Indeed, by the nonstandard characterization of the limit \cite{Dienerreeb} it holds that for all $ \varepsilon\gnsim 0$ there exists $ \delta\gnsim 0 $ such that $ \left| \frac{f(x)-f(a)}{x-a}-f'(a)\right| <\varepsilon $, while for some $ \eta\gnsim 0 $, whenever $ 0\lnsim |x-a|\leq\eta $, for some $ \theta $ with $ 0<\theta<1 $ it holds that $ \frac{f(x)-f(a)}{x-a}-f'(a)=\frac{x-a}{2}f''(a+\theta(x-a))\ncong 0 $. 
Let $ k\in \N,k\geq 2, a_{1}, \dots, a_{k}\in \R, a_{1},\dots, a_{k}\neq 0 $ be standard and $ P(x)=(a_{1}+\oslash) x+\cdots +(a_{k}+\oslash) x^{k} $. Then $ D_{\oslash}P(x)=\oslash + (2a_{2}+\oslash)x+\cdots +(ka_{k}+\oslash)x^{k-1}$ for limited $ x\in \R $ and $ D_{\oslash}P(x)=ka_{k}(1+\oslash)x^{k-1} $ for unlimited $ x\in \R $; then $ N(D_{\oslash}P(x))= \oslash x^{k-1}$, illustrating how the imprecision of $ P(x) $ increases with $ x $.

For $ M>0 $, the $ M $-derivative usually is an imprecise function, though the $ M $-derivative of a linear function reduces to the ordinary derivative. The following example concerns a flexible function $ F $ which is precise in one point, having still an imprecise $ M $-derivative. Indeed, let  $F(x)=x+x\cdot \oslash$ for all $x\in\R$ and $a\in \R.$ Then $ F(0)=0 $, while $\frac{dF}{d_{\oslash}x}(0)=1+\oslash$.

Consider the not everywhere continuous and differentiable function $ f:\R\rightarrow\R $ defined by
\begin{equation}\label{crep}
f(x)=\varepsilon\left( [x/\varepsilon]-x/\varepsilon\right)
\end{equation}
Note that it is neutrix-differentiable with respect ot non-zero neutrices, and at  sufficiently large scale the neutrix-derivative is neutricial, for instance $ \frac{df}{d_{\oslash}x}(0)=\lim\limits_{x\to \oslash}\frac{\varepsilon\left( [x/\varepsilon]-x/\varepsilon\right)}{x}=\varepsilon\pounds$ and $\frac{df}{d_{\varepsilon\pounds}x}(0)=\lim\limits_{x\to \varepsilon\pounds} \frac{\varepsilon\left( [x/\varepsilon]-x/\varepsilon\right)}{x}=\oslash$.

Next proposition states that the ordinary derivative is contained in the neutrix-derivative under a condition of inner continuity.
\begin{proposition}\label{ordneut}
Let $M$ be a neutrix,  $f : \ X\subseteq \R \rightarrow \E$ be an internal differentiable function and $a\in X$ be contained in an $M$-neighbourhood $ U\subseteq X $ of $a$. Assume that	 $ A\equiv {\rm{N}}\left( D_{M}f(a) \right)  $ is such that $f'$ is  $M\times A$-inner continuous at $a$. Then $ f'(a)\in D_{M}f(a) $.
\begin{proof}
 Because $f'$ is $M\times A$-inner continuous at $a$, for all $\eps>A$ there exists $\eta>M$ such that $ a+\eta \in U $ and $ \left|f'(a+k)-f'(a)\right|<\eps$ for all $ k $ with $ |k|<\eta $. It follows from the Mean Value Theorem that for all $ h $ with $M<|h|<\eta$ there exists $ \theta \in (0, 1)$ such that 

$$\left| \frac{f(a+h)-f(a)}{h}-f'(a)\right|=\left|f'(a+\theta h)-f'(a)\right|<\eps.$$
Hence 
$f'(a)\in A\mbox{-}\lim\limits_{h\to a+M} \frac{f(a+h)-f(a)}{h}=D_{M}f(a).$	
\end{proof}

\end{proposition} 
Theorem \ref{derop} indicates that, as in the case of limits, the $M$-derivative tends to satisfy the usual properties of algebraic operations with inclusions. 

\begin{theorem} \label{derop} Let $F, G$ be  flexible functions defined on $X\subseteq \R$ and $a\in X$ be an $M$-accumulation point of $X$.
	Assume that $F,G$ are $M$-differentiable at $a$. Then 
	\begin{enumerate}[(a)]
	\item \label{deradd} $ F\pm G$ is $M$-differentiable at $a$ and 	
	$$D_M(F\pm G)(a)\subseteq  D_MF(a)\pm  D_MG(a).$$	
	\item \label{dersca} If $ c\in \R $, $ c\in \E $ zeroless or $ D_MF(a)$ zeroless, the flexible function $cF $ is $M$-differentiable at $a$ and $$D_M (cF)(a)=c\cdot D_MF(a).$$ 
	\item \label{derprod} If $ F,G $ are $ M $-outer continuous at $ a $ and $ F(a) ,G(a)$  are zeroless, the flexible function $ FG$ is $M$-differentiable at $a$ and 	
	$$D_M(FG)(a)\subseteq  D_MF(a)G(a)+F(a)  D_MG(a).$$	
	\item \label{derdiv} If $ G $ is $ M $-outer continuous at $ a $ and $ G(a)$  is zeroless the flexible function $ 1/G$ is $M$-differentiable at $a$ and 	
	$$D_M\left( \frac{1}{G }\right)(a)\subseteq - \frac{D_MG(a)}{G^{2}(a)}.$$	
	\end{enumerate}	
	
\end{theorem} 
\begin{proof}
The properties follow in a straightforward way from Theorem \ref{operationsNconv}.
In Part \eqref{derprod} we apply the inclusion
\begin{align*}
\frac{F(a+h)G(a+h)-F(a)G(a)}{h}
\subseteq&\frac{F(a+h)-F(a)}{h}g(a+h)+\frac{{\rm{N}}_{F}(a)}{h}G(a)+\frac{{\rm{N}}_{F}(a)}{h}G(a+h)\\
&+f(a)\frac{G(a+h)-G(a)}{h}+F(a)\frac{{\rm{N}}_{G}(a+h)}{h}+F(a+h)\frac{{\rm{N}}_{G}(a+h)}{h},
\end{align*}
and, using \eqref{decdifM} and the $ M $-outer continuity of $ F $ and $ G $, we see that the limit for $ h\rightarrow M $ is included in $ D_MF(a)G(a)+F(a)  D_MG(a) $. In the proof of Part \eqref{derdiv} we apply the inclusion
\begin{equation*}
	\frac{1}{h}\left( \frac{1}{G(a+h)}-\frac{1}{G(a)}\right) 
	\subseteq-\frac{1}{h}\left(\frac{g(a+h)-g(a)}{g(a+h)g(a)}\right)
	+\frac{1}{G^{2}(a)}\frac{{\rm{N}}_{G}(a)}{h}+\frac{1}{G^{2}(a+h)}\frac{{\rm{N}}_{G}(a+h)}{h},
\end{equation*}
and again with the help \eqref{decdifM} and the $ M $-outer continuity of $ G $ we see that the limit for $ h\rightarrow M $ is included in $  D_MG(a)/G^{2}(a) $.
\end{proof}

We end this section with a Chain Rule.

\begin{theorem}\label{Tcr1}
	Let $ a\in \R$, $ M,N $ be neutrices, $\varphi: \R\rightarrow \R$ be $M$-differentiable at $a$ and $F: X\rightarrow \R$. Assume that $\varphi(a)$ is an $ N $-interior point of $X$, $\lim\limits_{x\to a+M}\varphi (x)=\varphi(a)+N$ and $ \varphi(x)\notin \varphi(a)+N $ in some $ M $-outer neighborhood $ U\subseteq X $ of $ a $, and $F$ is $N$-differentiable at $\varphi(a)$. If $ D_M\varphi(a) $ or $ D_N F(\varphi(a) )$ is zeroless, then $D_MF(\varphi(a))\subseteq D_N F(b )D_M\varphi(a)$.  
\end{theorem}
\begin{proof}  Because $\varphi$ is real-valued, $ F\circ\varphi $ is well-defined for $ x\in X$, and for $ x\in U\setminus a+M $
	\begin{equation*}
		\frac{F(\varphi(x))-F(\varphi(a))}{x-a}=\frac{F(\varphi(x))-F(\varphi(a))}{\varphi(x)-\varphi(a)}\frac{\varphi(x)-\varphi(a)}{x-a}.
	\end{equation*}
	By Theorem \ref{xylimits} $\lim\limits_{x\to a+M}\frac{F(\varphi(x))-F(\varphi(a)))}{\varphi(x)-\varphi(a)} 	\subseteq  \lim\limits_{y\to b +N}\frac{F(y)-F(b )}{y-b } $ is well-defined. Then by  Theorem \ref{operationsNconv}.\eqref{Theorem product}(i)
	\begin{align*}
		D_MF(\varphi(a))=&\lim\limits_{x\to a+M}\frac{F(\varphi(x))-F(\varphi(a))}{x-a}\subseteq \lim\limits_{x\to a+M}\frac{F(\varphi(x))-F(\varphi(a)))}{\varphi(x)-\varphi(a)}\lim\limits_{x\to a+M}\frac{\varphi(x)-\varphi(a)}{x-a}\\
		\subseteq & \lim\limits_{y\to b +N}\frac{F(y)-F(b )}{y-b }\lim\limits_{x\to a+M}\frac{\varphi(x)-\varphi(a)}{x-a}
		=  D_N F(b )D_M\varphi(a).
	\end{align*}
\end{proof}
\section{Nearly optimal points and neutrix derivatives}\label{section2c8}

By the classical Fermat Lemma the derivative of a differentiable function vanishes at an extreme point. If the involved neutrices are big enough, an $M$-differentiable function of a flexible function becomes neutricial at an $ L $-near-optimal point, giving a necessary condition for existence. Theorem \ref{necond} is formulated for near-minimizers.
\begin{theorem}\label{necond} Let $M$ be a neutrix and $X\subseteq \R, X\not=\emptyset$ and $ a $ be an $ M $-interior point of $ X $.
	Assume that $ F:X\rightarrow \R $ is $ M $-differentiable at $ a $. Let $ L $  be a neutrix which is stable for $ M $, and such that $L\supseteq \max\left( {\rm{N}}\left( F(a)\right), {\rm{N}}\left( \frac{dF}{d_M x}(a)\right)\right)  $.  If $a$ is an $M$-local $L$-minimizer of $F$,
	
	\begin{equation}
		\frac{dF}{d_M x}(a)\subseteq L.
	\end{equation}
	Moreover $ \frac{dF}{d_M x}(a)= L $ for $L={\rm{N}}\left( \frac{dF}{d_M x}(a)\right) $.
\end{theorem}

\begin{proof} Put  $K={\rm{N}}\left( \frac{dF}{d_M x}(a)\right) $. Then $ K\leq L $ and $\frac{dF}{d_M x}(a)=b+K$, with $ b\in \R $. Hence for all $\eps > K$ there exists $\delta_1>M$ such that for all $x\in X, M < |x-a|<\delta_1$ one has 
	\begin{equation} \label{eq4} \left| \frac{F(x)-F(a)}{x-a}-b+K\right| <\eps.
	\end{equation} 
	This implies that 
	\begin{equation} \label{eq6non} b+K-\eps<\frac{F(x)-F(a)}{x-a}+K
	\end{equation} 
	for  $M < a-x< \delta_1$ and 
	\begin{equation*}  \frac{F(x)-F(a)}{x-a}+K < b+K+\eps
	\end{equation*} 
	for  $M < x-a< \delta_1$, hence also 
	\begin{equation} \label{eq6nonthem} b+K+\eps> \frac{F(x)-F(a)}{x-a}+K. 
	\end{equation} 
	
	On the other hand, the point  $a$ is an $M$-local $L$-minimizer of $F$ on $X$, so there exists $ \delta_{2}>M $ such that $ F(x)\geq L+F(a) $ for all $x\in X, |x-a|\leq \delta_2$; then also 
	\begin{equation}F(x)- F(a)\geq L \label{eq1}, 
	\end{equation}
	because $ N(F(a))\subseteq L $. Put $\delta=\min\{\delta_1, \delta_2\}$. Note that $\oslash_L\subseteq M$, so for $|x-a|>M$ it holds that $1/(x-a)$ is not an exploder of $L$. For $x\in X, M < a-x<\delta$, dividing by $ x-a $ changes signs, so 
	
	$$\frac{F(x)- F(a)}{x-a}\leq \frac{L}{x-a}\leq L.$$  This implies that
	\begin{equation} \label{eq5} \frac{F(x)- F(a)}{x-a}+K\leq K+L=L
	\end{equation} 
	for all $x\in X, M < a-x<\delta$. It follows from \eqref{eq6non} and \eqref{eq5} that 
	\begin{equation} \label{eq7}
		b+K-\eps \leq L.
	\end{equation}	
	Similarly, we derive from \eqref{eq1} that $\frac{F(x)- F(a)}{x-a}\geq \frac{L}{x-a}\geq L$ for all $x\in X, M < x-a<\delta$. Then by \eqref{eq6nonthem} 
	\begin{equation} \label{eq8non}
		b+K+\eps\geq L
	\end{equation} 
	for all $x\in X, M < a-x<\delta$.
	
	Because $ \eps>K$ is arbitrary, formulas \eqref{eq7} and \eqref{eq8non} imply that $b+K\subseteq L$.   Indeed, if $b<L$, we choose $\eps=-b/2>L$, and then $b+K =b/2+K+\eps< L$, which is a contradiction to \eqref{eq8non}. If $b>L$, taking $\eps=b/2>L$ we obtain that $b+K-\eps=b/2+K>L$, which is contradictory to \eqref{eq7}. Hence $ b\in L $, and we conclude that 
	\begin{equation}\label{dFL}
		\frac{dF}{d_M x}(a) =b+ K\subseteq L+K=L.
	\end{equation}
	
	Finally, if $ {\rm{N}}\left( \frac{dF}{d_M x}(a)\right)=L$, we derive from \eqref{dFL} that $ \frac{dF}{d_M x}(a) \subseteq L= {\rm{N}}\left( \frac{dF}{d_M x}(a)\right)$. Hence $ \frac{dF}{d_M x}(a) = {\rm{N}}\left( \frac{dF}{d_M x}(a)\right)= L$.
\end{proof}

Theorem \ref{necond} permits also to recognize points which do not correspond to extremes of a flexible function $ F $. Indeed, if an $M$-derivative of $F$ is zeroless at some point $ x_{1} $, this point is not an $M$-local $K$-extreme point of $F$, for $ K={\rm{N}}\left( D_{M}F(a)\right)  $. 
\begin{examples}Let $F : \ \R \longrightarrow \E$ be given by  $F(x)=x^2+\oslash$ for all $x\in \R$. Then $\frac{d_{\oslash} F}{d_{\oslash} x}(x)=2x+\oslash.$ So  $\frac{d_{\oslash} F}{d_{\oslash} x}(x)
   	=2x+\oslash \not \subseteq \oslash$ for all $x\in \R,x\not\in \oslash$, hence these points are not $\oslash$-local $\oslash$-minimizers of $F$.  At the obvious  $\oslash$-minimizer $ 0 $, the neutrix $ L=\oslash $ satisfies $\oslash\supseteq {\rm{N}}\left( \frac{dF}{d_M x}(0)\right)=\oslash $, $\oslash\supseteq {\rm{N}}\left( F(0)\right)=\oslash $ and $ \oslash $ also contains the absorbers of $ M=\oslash $, and we have indeed $ \frac{d F}{d_{\oslash} x}(0)=\frac{d_{\oslash} F}{d_{\oslash} x}(0)
   	= \oslash$; also infinitesimals are $\oslash$-minimizers, and the above relations remain true for $ x\simeq 0 $.
   	
   	Define $g : \ \R \longrightarrow \E$ by $ g(x)=x^{2} + f(x) $, with $ f $ given by \eqref{crep}. Then both $ \frac{dg}{d_{\oslash}x}(0)=2x+\oslash+\varepsilon\pounds=2x+\oslash$ and $\frac{dg}{d_{\varepsilon\pounds}x}(0)=2x+\varepsilon\pounds+\oslash=2x+\oslash$, and Theorem \ref{necond} is satisfied with $ L=\oslash $ for respectively $ x\in \oslash $ and  $x\in \varepsilon\pounds $ with equality.
\end{examples}

\section{Neutrix-partial derivatives }\label{section9c6}
\begin{remark}
	In the remaining sections we study flexible functions of several variables from the point of view of differentiability and optimization. Though some results are easily generalized to an arbitrary standard number of variables, for reasons of simplicity we consider the case of two variables.
\end{remark}
We define partial derivatives by extending Definition \ref{defderNM}. We introduce also a notion of \emph{strong differentiability}, giving rise to a kind of differential. The latter is needed to formulate a chain rule.

\begin{definition}  Let  $n\in \N$ be standard, $X\subseteq \R^2$ and $F :  \ X \longrightarrow \ \E$ be a flexible function defined on $X$. Let $M_i, N_i$ be neutrices for $i\in\{1,2\}$ and let $(a, b)\in \R^2$  be a $(M_{1}\times M_{2})  $-accumulation point of $ X $. For $ i\in\{1,2\} $	the $M_i\times N_i$-partial derivative $\frac{\partial_{N_i} F}{\partial_{M_i} x}$ at $ (a,b) $ is defined by
	\begin{align*}
		\frac{\partial_{N_1} F}{\partial_{M_1} x}(a, b)&=N_1\mbox{-}\lim\limits_{x\to a+M_1}\frac{F(x, b)-F(a, b)}{x-a}\\
		\dfrac{\partial_{N_2} F}{\partial_{M_2} y}(a, b)&=N_2\mbox{-}\lim\limits_{y\to b+M_2}\dfrac{F(a, y)-F(a, b)}{y-a}.
	\end{align*}
	In analogy to functions of one variable, taking $ N_{1} $ respectively $ N_{2} $ minimal we get  the notions of $M_1$-partial derivative respectively $M_2$-partial derivative as the outer limits
	\begin{align*}
		\dfrac{\partial F}{\partial_{M_1} x}(a, b)&=\lim\limits_{x\to a+M_1}\dfrac{F(x, b)-F(a, b)}{x-a}\\
		\dfrac{\partial F}{\partial_{M_2} y}(a, b)&=\lim\limits_{y\to b+M_2}\dfrac{F(a, y)-F(a, b)}{y-a}.
	\end{align*}
\end{definition}
Theorem \ref{ferma ham nhieu bien} extends Theorem \ref{necond} to flexible functions of two variables.

\begin{theorem}\label{ferma ham nhieu bien} 
	Let $M=(M_1, M_2)$ be a neutrix. Let $ X\subseteq \R^2 $, $a=(a_1, a_2)\in X$ be an $M$-interior point of $X$ and $F :  X  \rightarrow \E$ be a flexible function which is $M_1$-partial differentiable in the first variable and $M_2$-partial differentiable in the second variable at $a$. Let $L$ be a neutrix which is stable for $\min\{M_1, M_2\}$, such that $ L\supseteq \max\left( N(F(a)),{\rm{N}}\left(\frac{\partial F}{\partial_{M_1} x_1}(a) \right),{\rm{N}}\left(\frac{\partial F}{\partial_{M_2} x_2}(a) \right)\right) $.	If   $a$ is an $M$-local $L$-extreme point of $F$, then $\frac{\partial F}{\partial_{M_i} x_i}(a)\subseteq L$ for $i\in\{1,2\} $.
\end{theorem}
\begin{proof}
	Put $G_{1}(x_{1})=F(x_1, a_{2})$, $G(x_2)=F(a_1, x_{2})$. Since $a$ is an $M$-local $L$-extreme point of $F$, for $i\in\{1,2\}$ it holds that $a_i$ is an $M_i$-local $L$-extreme point of $G_{i}$, so $\frac{\partial F}{\partial_{M_i} x_i}(a)=\frac{d G_{i}}{d_{M_i} x_i}(a_i)\subseteq L$ by Theorem \ref{necond}. 
\end{proof}

In classical analysis a function $f(x, y)$ is differentiable at $(a_1, a_2)\in \R^{2}$ if their exist "error-functions" $ \alpha_1,\alpha_2 $ such that $f(a_1+h_1, a_2+h_2 )-f(a_1, a_2)=f_x(a_1, a_2)h_1+f_{y}(a_1, a_2)h_2+\alpha_1(h_1, h_2)h_1+\alpha_2(h_1, h_2)h_2$ on some neighbourhood of $(a_1 a_2)$, where $\lim\limits_{(h_1, h_2)\to (0,0)}\alpha_1(h_1, h_2)=\lim\limits_{(h_1, h_2)\to (0, 0)}\alpha_2(h_1, h_2)=0$. In Definition \ref{totdif} we extend this notion to flexible functions using the mixed limit of Definition \ref{milim}. Indeed, consider the external point $ (a_{1}+M_{1} ,a_{2}+M_{2})$, where $ M_{1} ,M_{2}$ are neutrices. To define the neutrix-partial derivative in the first variable, this variable should stay outside $ a_{1}+M_{1} $, while it is useful to have information on the error-functions $ \alpha_1(x,y), \alpha_2(x,y) $  for $ y $ inside $  a_{2}+M_{2}$, and vice-versa. In this way we will be able to extend the Chain Rule of Theorem \ref{Tcr1} to functions of two variables, see Theorem \ref{chain} below.

\begin{definition}\label{totdif} Let $M=(M_1, M_2)$ be a neutrix. 
	Let $F:X\ \to \E$ where $X\subseteq \R^2$ be a flexible function and $a=(a_1, a_2)$ be an $M$-interior point of $X$.  Assume that $F$ has $M_i$-partial derivatives at $a$ for $i\in\{1,2\} $, and on some $M$-neighborhood $U$ of $a$ flexible functions $ \alpha_{1},\alpha_{2} $ are given such that for all $h=(h_1, h_2)\neq (0, 0)$ such that  $a+h\in U$
	\begin{equation}\label{defdiff} 
		F(a_1+h_1, a_2+h_2)-F(a_1, a_2)-\dfrac{\partial F}{\partial_{M_1} x_1}(a_1, a_2) h_1-\dfrac{\partial F}{\partial_{M_2} x_2} (a_1, a_2)h_2=\alpha_1(h_1, h_2)h_1+\alpha_2(h_1 h_2) h_2.
	\end{equation} 
	\begin{enumerate}[(a)]
		\item The flexible function $F$ is said to be \emph{strongly $M$-differentiable} with respect to $x_1$ at $a$ if 
		$$\lim\limits_{\substack{h_1\to M_{1}\\ h_2\twoheadrightarrow M_2}}\alpha_1(h_1, h_2)=N\left(\dfrac{\partial F}{\partial_{M_1} x_1}\right) ,\lim\limits_{\substack{h_1\to M_{1}\\ h_2\twoheadrightarrow M_2}}\alpha_2(h_1, h_2)= N\left(\dfrac{\partial F}{\partial_{M_2} x_2}\right).$$	
		
		\item 	The flexible function $F$ is said to be \emph{strongly $M$-differentiable} with respect to $x_2$ at $a$ if  	 
		$$\lim\limits_{\substack{h_1\twoheadrightarrow M_{1}\\ h_2\to M_2}}\alpha_1(h_1, h_2)=N\left(\dfrac{\partial F}{\partial_{M_1} x_1}\right) ,\lim\limits_{\substack{h_1\twoheadrightarrow M_{1}\\ h_2\to M_2}}\alpha_2(h_1, h_2)= N\left(\dfrac{\partial F}{\partial_{M_2} x_2}\right).$$			
	\end{enumerate}
\end{definition}

\begin{remark} 
	The equality \eqref{defdiff} becomes an inclusion if we put the partial derivatives at the right-hand side, and we have for all $h=(h_1, h_2)\neq (0, 0) $ such that $ a+h\in U, $
	\begin{equation*}
		F(a_1+h_1, a_2+h_2)-F(a_1, a_2) \subseteq \dfrac{\partial F}{\partial_{M_1} x_1}  h_1+ \alpha_1(h_1, h_2)h_1 + \dfrac{\partial F}{\partial_{M_1} x_2}  h_2+ \alpha_2(h_1, h_2)h_2 .
	\end{equation*}	
\end{remark}

\begin{example}\label{twovariables} Let $ F : \R^{2}\rightarrow\R $ be defined by $F(x, y)=x+y^2+\oslash$. Then $ F $ is strongly $(\oslash, \oslash)$-differentiable with respect to $x$ at $(1, 1)$. Indeed, $\frac{\partial F}{\partial_\oslash x}(1, 1)=1+\oslash$ and $\frac{\partial F}{\partial_\oslash y}(1, 1)=2+\oslash$, and
	\begin{align*} F(1+h_1, 1+h_2)-F(1, 1)&-\dfrac{\partial F}{\partial_\oslash x}(1, 1)-\dfrac{\partial F}{\partial_\oslash y}(1, 1)
		=h_1\oslash +h_2^2 +\oslash h_2+\oslash\\
		&=(\oslash+\oslash/h_1) h_1+(h_2+\oslash)h_2=\alpha_1(h_1, h_2)h_1+\alpha_2(h_1, h_2)h_2. 
	\end{align*}  
	with  $\alpha_1(h_1, h_2)=\oslash+\oslash/h_1$ for $ h_{1}\neq 0 $ and $\alpha_2(h_1, h_2)=h_2+\oslash. $ 	
	Then $$\lim\limits_{\substack{h_1\to \oslash\\ h_2\twoheadrightarrow \oslash}}\alpha_1(h_1, h_2)=\oslash=N\left(\dfrac{\partial F}{\partial_\oslash x}(1, 1)\right), \hspace{0,5cm}\lim\limits_{\substack{h_1\to \oslash\\ h_2\twoheadrightarrow \oslash}}\alpha_2(h_1, h_2)=\oslash =N\left(\dfrac{\partial F}{\partial_\oslash y}(1, 1) \right).$$
	Similarly, using $\beta_1(h_1, h_2)=\oslash$ and  $\beta_2(h_1, h_2)=h_2+\oslash+\oslash/h_2 $ for $ h_{2}\neq 0 $ we see that $ F $ is strongly $(\oslash, \oslash)$-differentiable with respect to $y$ at $(1, 1)$.	
\end{example}

\begin{theorem}\label{chain}
	 Let  $ X\subseteq\R $, $\varphi\equiv(y_{1},y_{2}): X\rightarrow \R^2$, let $F:\varphi(X)\rightarrow \E$ be a flexible function and let $G=F\circ  \varphi$.  Let $ M $ be a neutrix, $a$ be an $M$-interior point of $X$ and $ b=(b_1, b_2) \equiv \left(  y_1(a), y_2(a)\right) $. Assume that 
	\begin{enumerate}[(a)]

		\item \label{condlimit}	There exist neutrices $ N_{1},N_{2} $ such that $ b $ is an $ N $-interior point of $ \varphi(X) $, where $N\equiv (N_{1}, N_{2})$, and for  $ i\in\{1,2\} $ 
		 \begin{align}\label{Nilim}
			N_{i}\mbox{-}\lim_{h\twoheadrightarrow M} y_i(a+h)= b_i+N_{i}.
		\end{align}
	\item  There exists $p\in\{1,2\}$ and an $M$-neighbourhood $U$ of $a+M$ such that whenever $a+h\in U, h\notin M$ 	
			\begin{equation}\label{outside}   
			y_p(a+h)\notin b_p+N_{p},
		\end{equation}
		and $F$ is strongly $N$-differentiable with respect to $y_p$ at $b$.

\item\label{ypd}   $\frac{dy_1}{d_M x}(a), \frac{dy_2}{d_M x}(a), \frac{\partial F}{\partial_{N_1} y_1}(b_1, b_2)$ and $\frac{\partial F}{\partial_{N_2} y_2}( b_1, b_2)$ are well-defined, and $\frac{dy_1}{d_M x}(a)$ or $\frac{\partial F}{\partial_{N_1} y_1}(b_1, b_2)$ is zeroless,  and also  $\frac{dy_2}{d_M x}(a)$ or $\frac{\partial F}{\partial_{N_2} y_2}( b_1, b_2)$ is zeroless.  
	\end{enumerate}
	Then 
	\begin{align*}
		\frac{d G}{d_{M} x}(a)\subseteq &\frac{\partial F}{\partial_{N_{1}} y_1}(b_1,b_2 )  \frac{d y_1}{d_{M} x}(a) +\frac{\partial F}{\partial_{N_{2}} y_2}(b_1,b_2)  \frac{d y_2}{d_{M} x}(a).
	\end{align*} 
	
\end{theorem}

\begin{proof} Without loss of generality, in assumption \eqref{outside} we assume $p=1$. Let $ k=(k_{1},k_{2}) $, where for $ i\in\{1,2\} $
	$$k_i\equiv k_{i}(h)=y_i(a+h)-y_i(a),$$	
	otherwise said $y_i(a+h)=b_i+k_i$. Because $ y_{1},y_{2} $ are $ M $-differentiable at $a$, for $ i\in\{1,2\} $
	\begin{equation}\label{partialy}
		\lim_{h\to M}	\frac{k_i(h)}{h}= \frac{d y_i}{d_M x}(a). 
	\end{equation}
	Because $F$ is strongly $ (N_{1}, N_{2})$-differentiable w.r.t $y_1$ at $(b_1, b_2)$, there exist flexible functions $\alpha_1, \alpha_2$ defined on some $ N$-neighborhood of $ (b_1, b_2) $ such that
	$\lim\limits_{\substack{u_1\to N_1\\u_2\twoheadrightarrow N_2}} \alpha_1(u_1, u_2)=N\left( \frac{\partial F}{\partial_{N_{1}} y_1}(b_1,b_2 )\right) \equiv K_{1}$ and $\lim\limits_{\substack{u_1\to N_1\\u_2\twoheadrightarrow N_2}}\alpha_2(u_1, u_2)=N\left( \frac{\partial F}{\partial_{N_2} y_2}(b_1,b_2 )\right) \equiv K_{2}$ and 

	\begin{align}\label{cm1}
	G(a+h)-G(a)=&F\big(y_1(a+h), y_2(a+h)\big)-F(y_1(a), y_2(a)) \nonumber 
	=	F(b_1+k_1, b_2+k_2)-F(b_1, b_2)\\
	\subseteq &\frac{\partial F}{\partial_{N_1} y_1}(b_1,b_2 ) k_1+\frac{\partial F}{\partial_{N_2} y_2}(b_1, b_2) k_2+\alpha_1(k_1, k_2)  k_1+\alpha_2(k_1, k_2) k_2 \notag\\
	=&\left( \frac{\partial F}{\partial_{N_1} y_1}(b_1,b_2 )+\alpha_1(k_1, k_2) \right)  k_1+\left( \frac{\partial F}{\partial_{N_2} y_2}(b_1, b_2) +\alpha_2(k_1, k_2) \right) k_2.
\end{align}	
By Theorem \ref{xylimits} and Remark \ref{oim}
\begin{equation}\label{inclusion neutrix 1}
	\lim\limits_{h\to M}\alpha_1(k_1, k_2) =	\lim\limits_{h\to M}\alpha_1(k_1(h), k_2(h))\subseteq \lim\limits_{\substack{u_1\to N_1\\u_2\twoheadrightarrow N_2}} \alpha_1(u_1, u_2) =K_1
\end{equation}
\begin{equation}
	\label{inclusion neutrix 2}\lim\limits_{h\to M}\alpha_2(k_1, k_2) =	\lim\limits_{h\to M}\alpha_2(k_1(h), k_2(h))\subseteq \lim\limits_{\substack{u_1\to N_1\\u_2\twoheadrightarrow N_2}} \alpha_2(u_1, u_2) =K_2.
\end{equation}	
It follows from  assumption \eqref{ypd}, \eqref{cm1}, Theorem \ref{operationsNconv},  \eqref{partialy}, \eqref{inclusion neutrix 1}, and \eqref{inclusion neutrix 2} that 
\begin{align*}\label{cm2} 
	\frac{\partial G}{\partial_{M} x}(a )=&	\lim_{h\to M} \frac{G(a+h)-G(a)}{h}\\
	\subseteq & \lim_{h\to M}\frac{\left( \frac{\partial F}{\partial_{N_1} y_1}(b_1,b_2 )+\alpha_1(k_1, k_2) \right)  k_1+\left( \frac{\partial F}{\partial_{N_2} y_2}(b_1, b_2) k_2+\alpha_2(k_1, k_2) \right) k_2}{h} \\
	\subseteq &\lim\limits_{h\to M} \left( \frac{\partial F}{\partial_{N_1} y_1}(b_1,b_2 )+\alpha_1(k_1, k_2) \right) \lim_{h\to M} \frac{k_1}{h}+ \frac{\partial F}{\partial_{N_2} y_2}(b_1, b_2)  \lim_{h\to M}\frac{k_2}{h}\\
	\subseteq &\left(\frac{\partial F}{\partial_{N_1} y_1}(b_1,b_2 ) +K_1\right) \frac{d y_1}{d_M x}(a) + \left(\frac{\partial F}{\partial_{N_2} y_2}(b_1,b_2) +K_2\right)  \frac{d y_2}{d_M x}(a)		 \\
	=& \frac{\partial F}{\partial_{N_1} y_1}(b_1,b_2 )  \frac{d y_1}{d_M x}(a) +\frac{\partial F}{\partial_{N_2} y_2}(b_1,b_2)  \frac{d y_2}{d_M x}(a).  
\end{align*}

\end{proof}

   \section{An implicit function theorem for  neutrix-differentiability}\label{section11c6} 
Let	$M, N$ be neutrices and $ V $ be an $(M, N)$-neighborhood of a point $(a, b)\in \R^{2}$. Let $g$ be an internal real function defined on $V$, which is $ N $-differentiable with respect to the second variable. We give conditions such that an $M$-differentiable implicit function exists on some $M$-neighborhood of $a $. The conditions include ordinary differentiability of $ g $, but for the rest are essentially imprecise. For instance $ (a,b) $ should not be a "nearly singular" point of $ g $, in the sense that 
\begin{equation}\label{defgamma}
\gamma_x\equiv \lim\limits_{\substack{(x,y)\twoheadrightarrow (a+M,b+N)}}\frac{\partial g}{\partial x}(x, y), \hspace{1cm} \gamma_y\equiv \lim\limits_{\substack{(x,y)\twoheadrightarrow (a+M,b+N)}}\frac{\partial g}{\partial y}(x, y)
\end{equation} are zeroless. The implicit function will be used in Section \ref{seclag} to prove the existence of Lagrange multipliers in a near-minimizing problem.

\begin{theorem}\label{implicitfunction2}
Let $(a, b)\in \R^2$, $M, N$ be neutrices and $V$ be an $(M, N)$-neighborhood of $(a, b)$. Let	$g:V\rightarrow \R$ be an internal differentiable function such that $g(a, b)=0$. Consider the equation $g(x, y)=0$ for $ (x,y)\in V $. Let $\gamma_x,\gamma_y  $ be given by \eqref{defgamma}. Assume that $ \frac{\partial g}{\partial_N y}(a, b)$ is well-defined and zeroless, where $ A\equiv {\rm{N}}\left( \frac{\partial g}{\partial_N y}(a, b)\right)  $ is such that 

\begin{enumerate}[(a)]
	\item \label{zerolessassumption2} $ N$ is stable for $A$.	
	\item $g, \frac{\partial g}{\partial y}$ are  $(M, N)\times A$-inner continuous on $V$.
	\item $ \gamma_y$ is zeroless.
\end{enumerate}
Then  there exist $ \delta>M $ and a unique $ M $-differentiable function $f: U \equiv[a-\delta, a+\delta]\rightarrow \R$ such that $g(x,f(x))=0$ for  $x\in U$, $f(a)=b$,  and
	\begin{equation}\label{DMf}
		D_Mf(a)\subseteq -\frac{\gamma_x}{\gamma_y}.
	\end{equation}
In addition, if $\gamma_x$  zeroless and $\left( \gamma_x/ \gamma_y\right)  M\subseteq N$, the function $ f $ is $ M\times N $-inner continuous at $ a $.
\end{theorem}

\begin{proof} 

Firstly, we show that there exists and $ \delta>M , d>N $ such that for all $ x\in U \equiv[a-\delta, a+\delta] $
	\begin{equation}\label{ddelta}
	g(x, b-d)<A,g(x, b+d)>A,
	\end{equation}
which will imply the existence of the implicit function $ h $.
By Proposition \ref{ordneut} the ordinary derivative $\frac{\partial g}{\partial y}(a, b)$ is contained in the $ N $-derivative $ \frac{\partial g}{\partial_N y}(a, b)$. Since $\frac{\partial g(a, b)}{\partial_N b}=\frac{\partial g}{\partial y}(a, b)+A$ is zeroless, without loss of generality we may assume that  $\frac{\partial g(a, b)}{\partial_N b}>A$. Also  $\frac{\partial g}{\partial y}$ is $(M, N)\times A$-continuous at $(a, b)$, so there exist $c_1>M, c_2>N$ such that $\frac{\partial g}{\partial y}(x, y)>A$ for all $(x, y)\in [a-c_1, a+c_1]\times [b-c_2, b+c_2]$.  As a consequence, $\frac{\partial g}{\partial y}(x, y)>0$ for all $(x, y)\in [a-c_1, a+c_1]\times [b-\alpha_2, b+\alpha_2]$. In addition, $N$ contains all absorbers of $A$, so Theorem \ref{necond} and the fact that $g(a, b)=0$ imply that there exists $d>N$ such that 
	
	\begin{equation}\label{ginc}
		\begin{array}{ll} 
			g(a, b-d)&=g(a, b-d)-g(a, b)=g(a, b-d)=\varphi(b-d)<A   \\
			g(a, b+d)&=g(a, b+d)-g(a, b)=\varphi(b+d)>A,
		\end{array}   
	\end{equation}
	where we may assume that $d<\alpha_1, d<\alpha_2$. Because $g$ is $(M, N)\times A$-continuous on $V$, the functions $ \psi_1,\psi_2 $ given by $\psi_1(x)\equiv g(x, b-d)$ and $\psi_2(x)\equiv g(x, b+d)$ are $M\times A$-continuous at $a$.  Let $\eps_1=-g(a, b-d)/2$. By \eqref{ginc}, $\eps_1>A$. Since $g(x, b-d)$ is $M\times A$-continuous at $a$, there exists $\delta_1>M$ such that  $|g(x, b-d)-g(a, b-d)|<\eps_1$ for $|x-a|\leq \delta_1$. As a consequence, $g(x, b-d)<g(a,b-d)/2<A$ for $|x-a|\leq \delta_1$. Similarly, there exists $\delta_2>M$ such that $g(x, b+d)>A$ for all $|x-a|\leq \delta_2$. Let $\delta=\min\{\delta_1, \delta_2\}$, then \eqref{ddelta} holds for all $x\in [a-\delta, a+\delta]$. 
	
	As a consequence $g(x, b-d)<0, g(x, b+d)>0$ holds for all $x\in [a-\delta, a+\delta]$. 	Using the intermediate value theorem we conclude that there exists a unique function  $y=f(x)$ determined on $[a-\delta, a+\delta]$ with $\delta>M$ which is continuously differentiable on this interval and $f(a)=b$. 
	
	Secondly, we prove \eqref{DMf}. 	As long as $(a+h, f(a+h)) \in V $,
	\begin{align*}
		0= &	g\big(a+h, f(a+h))-g(a, f(a)\big)=g(a+h, f(a+h))-g(a, f(a+h))+g(a, f(a+h))-g(a, f(a))\\
		= & \frac{\partial g}{\partial x}(a+\theta_1h, f(a+h))h+\frac{\partial g}{\partial y}(a, f(a)+\theta_2(f(a+h)-f(a)))(f(a+h)-f(a)),
	\end{align*}	
	where $\theta_1=\theta_1(h), \theta_2=\theta_2(h)\in (0, 1)$; for simplicity of notation we will not write the functional dependence. Hence
	\begin{equation}\label{a represent of quotient}  \frac{f(a+h)-f(a)}{h}=-\frac{\frac{\partial g}{\partial x}\big(a+\theta_1h, f(a+h)\big)}{\frac{\partial g}{\partial y}\big(a, f(a)+\theta_2(f(a+h)-f(a))\big)}.
	\end{equation} 
	Observe that 
	\begin{equation}\label{limgamma}
		\lim\limits_{h\twoheadrightarrow  M}\frac{\partial g}{\partial x}(a+\theta_1 h, f(a+h))\subseteq \gamma_x,
		\hspace{1cm}
		\lim\limits_{h\twoheadrightarrow  M}\frac{\partial g}{\partial y}(a, f(a)+ \theta_2(f(a+h)-f(a)) \subseteq \gamma_y,
	\end{equation}	
Then it follows from \eqref{a represent of quotient}, Theorem \ref{operationsNconv}  and Theorem \ref{inclusion of both limits} that $D_Mf(a)=\lim\limits_{h\to M} \frac{f(a+h)-f(a)}{h}\subseteq  -\frac{\gamma_x}{\gamma_y}.$

	Finally we prove that $ f $ is $ M\times N$-inner continuous. By \eqref{a represent of quotient}
	 \begin{equation}\label{repf}
		f(a+h)-f(a)=-\frac{\frac{\partial g}{\partial x}\big(a+\theta_1h, f(a+h)\big)}{\frac{\partial g}{\partial y}\big(a, f(a)+\theta_2(f(a+h)-f(a))\big)}h.
	\end{equation}
	Because the limits given by \eqref{limgamma} are zeroless, it follows from Theorem \ref{operationsNconv} that 
\begin{align*}\label{inclusin for f a}
	\lim\limits_{h\twoheadrightarrow M}	f(a+h)-f(a)\subseteq-\frac{\lim\limits_{h\twoheadrightarrow  M}\frac{\partial g}{\partial x}(a+\theta_1 h, f(a+h))}{\lim\limits_{h\twoheadrightarrow  M}\frac{\partial g}{\partial y}\big(a, f(a)+\theta_2(f(a+h)-f(a))\big)}\lim\limits_{h\twoheadrightarrow  M}h
	\subseteq -\frac{\gamma_x}{\gamma_y}M\subseteq N.
\end{align*} 	
Hence $ f $ is $ M\times N $-inner continuous at $ a $.

\end{proof}
\begin{example}
	Let $g(x, y)=1-x^2-y^2$. Then $\frac{\partial g}{\partial x}(x, y)=-2x, \frac{\partial g}{\partial y}(x, y)=-2y,\gamma_x=\oslash, \gamma_y=-2+\oslash, $ and $g$ satisfies all assumptions of Theorem \ref{implicitfunction2} with $M=N=\oslash$ and $(a, b)=(0, 1)$. In fact $y=f(x)=\sqrt{1-x^2}$,  so $D_\oslash f(0)=\oslash$, and also	$-\frac{\gamma_x}{\gamma_y}=-\frac{\oslash}{2+\oslash}=\oslash$. 
\end{example}

\section{Lagrange multiplier}\label{seclag}

The Lagrange multiplier method for conventional optimization problem with an objective function $ f:\R^{n} \rightarrow \R$ and constraints $ g_{1}(x_1,x_2, \dots, x_n)=0,\dots, g_{m}(x_1,x_2, \dots, x_n)=0 $ asserts that if $a\in \R^n$ is a minimizer, there exist multipliers $\lambda_1, \dots, \lambda_m$ such that  
$$\begin{cases}\frac{\partial f}{\partial x_j}(a)-\sum\limits_{i=1}^m \lambda_i \frac{\partial g_i}{\partial x_j}(a)=0\hspace{1cm} (j=1, \dots, n)\\
	g_i(a)=0.\end{cases}$$
We consider the Lagrange multiplier method for a flexible function $F$ of two variables, with one internal constraint $ g=0 $, i.e. the optimization problem 
\begin{subequations}\label{lagrangetq} 
	\begin{equation}\label{2Lagrange problem} 
		\min_{(x, y)\in \R^2}F(x, y)
	\end{equation}  
	subject to the constraint 
	\begin{equation}\label{1Lagrange constraint}
		g(x,y)=0.
	\end{equation} 
\end{subequations}

The Lagrange multipliers are an indicator of changes in the value of the constraints. In our context we take also into account situations in which, due to all kind of uncertain circumstances, the objective is only known within some range of imprecisions, while also the effect of changes in the circumstances are only approximately known; this is modelled by partial derivatives of $ F $ with respect to neutrices $ M,N $. This means that though the constraint, say a budgetary restriction, is rigid, an optimum can only be approximate. We show that there exists a Lagrange multiplier $ \lambda $ and  neutrices $K_{1},K_{2}$ such that
\begin{equation*}		
	\begin{cases}\frac{\partial F}{\partial_{M} x}(a,b)-\lambda \frac{\partial g}{\partial x}(a, b)\subseteq  K_1\\
		\frac{\partial F}{\partial_{N} y}(a,b)-\lambda \frac{\partial g}{\partial y}(a, b)\subseteq  K_{2}
	\end{cases}.
\end{equation*}
One of the neutrices may be chosen to be neutrix of a partial derivative of $ F $, but for the other neutrix a correction term must be added, depending on the neutrix of the remaining partial derivative of $ F $, and the partial derivatives of $ g $. The proof uses the implicit function of Theorem \ref{implicitfunction2}, so the Lagrange multiplier is given outside nearly singular points of the constraints.

\begin{theorem}[Main Theorem]\label{Lagrange multiplier}  Let $M, N $ be neutrices. Consider the problem \eqref{2Lagrange problem}-\eqref{1Lagrange constraint}. Assume that the objective function $F$ is strongly  $(M,N)$-differentiable with respect to $x$ at $(a, b)$, the constraint $g$ satisfies the assumptions of Theorem \ref{implicitfunction2} and $(a, b)$ is an $(M, N)$-local $L$-minimizer of the problem \eqref{2Lagrange problem}-\eqref{1Lagrange constraint}, with 
	\begin{equation}\label{defL}
		L\equiv {\rm{N}}\left(\frac{\partial F}{\partial_{M} x}(a, b)\right) +  {\rm{N}}\left( \frac{\gamma_x}{\gamma_y}\frac{\partial F}{\partial_{N} y}(a, b)\right).
	\end{equation}
If $ L $ is stable for $ M  $ and contains $  {\rm{N} }(F(a, b))$, there exists   $\lambda\in \R$ such that 
		
		\begin{equation}\label{lagid}		
		\begin{cases}\frac{\partial F}{\partial_{M} x}(a,b)-\lambda \frac{\partial g}{\partial x}(a, b)\subseteq  L\\
		\frac{\partial F}{\partial_{N} y}(a,b)-\lambda \frac{\partial g}{\partial y}(a, b)= {\rm{N}}\left( \frac{\partial F}{\partial_{N}y}(a, b)\right)
	\end{cases}.
\end{equation}

\end{theorem}
\begin{proof}
	Let
	$\lambda\in   \frac{\partial F}{\partial_{N} y}(a, b)/\frac{\partial g}{\partial y}(a, b), $ which is well-defined because $\frac{\partial g}{\partial y}(a, b)\neq 0$. Then  $\lambda \frac{\partial g}{\partial y}(a, b)\in \frac{\partial F}{\partial_{N} y}(a, b)$, or
	\begin{equation}\label{ptn1}  \frac{\partial F}{\partial_{N} y}(a, b)-\lambda \frac{\partial g}{\partial y}(a, b) ={\rm{N}}\left( \frac{\partial F}{\partial_{N}y}(a, b)\right), 
	\end{equation} 
	which implies the second part of \eqref{lagid}. 
	
	To prove the first part, 
by Theorem \ref{implicitfunction2} the constraint $g(x, y)=0$ determines an implicit function $y=h(x)$, with well-defined $M$-derivative at $a$. Using the the notations of \eqref{defgamma} and Theorem \ref{inclusion of both limits}\eqref{limin} we get
	
	\begin{equation}\label{derh}
	\frac{d  h}{d_{M} x} (a)\subseteq  -\gamma_x/\gamma_y =-\frac{\partial g}{\partial x}(a, b)\Big{\slash} \frac{\partial g}{\partial y}(a, b)+{\rm{N}}\left(\gamma_x/\gamma_y\right),
	\end{equation}
which is zeroless. We will show that 
\begin{align}\label{ptn2} \frac{\partial F}{\partial_{M} x}(a, b) - \lambda \frac{\partial g}{\partial x}(a, b)	
	\subseteq \frac{\partial  F}{\partial_{M} x}(a, b)+ \frac{\partial  F}{\partial_{N} y}(a, b)\frac{d h}{d_{M} x}(a)- \lambda\left(  \frac{\partial g}{\partial x}(a, b)
	+  \frac{\partial g}{\partial y}(a, b) \frac{d h}{d_{M} x}(a)\right) ,
\end{align} 
with both 
	\begin{equation}\label{FL}
		\frac{\partial  F}{\partial_{M} x} (a, b)+  \frac{\partial  F}{\partial_{N} y}(a, b)\frac{d  h}{d_{M} x} (a)\subseteq L
	\end{equation}
and \begin{equation}\label{gL}
\lambda\left(  \frac{\partial g}{\partial x}(a, b)
+  \frac{\partial g}{\partial y}(a, b) \frac{d h}{d_{M} x}(a)\right)\subseteq L.
\end{equation}

We prove first \eqref{FL}.   Let $y_1(x)=x$ and $y_2(x)=h(x)$. Then $y_1(x), y_2(x)$ satisfy the assumptions of Theorem \ref{chain} with $N_1=M$, $N_2=N$ and $p=1$. Indeed, firstly $M\mbox{-}\lim\limits_{x\twoheadrightarrow a+M}y_1(x)=a+M$ and by Theorem \ref{implicitfunction2}, $h$ is $M\times N$-inner continuous at $a$ and $h(a)=b$, hence $N\mbox{-}\lim\limits_{x\twoheadrightarrow a+M} y_2(x)=N\mbox{-}\lim\limits_{x\twoheadrightarrow a+M}h(x)=b+N$. Secondly, $y_1(a+h)-y_1(a)=h\notin M$ if $h\notin M$ and $F$ is strongly $(M, N)$-differentiable with respect to $x$ at $(a, b)$ by assumption. Finally $\frac{dy_1}{d_M x}(a)=1$ is zeroless, and also $\frac{d y_2}{d_M x}(a) =\frac{dh}{d_M x}(a)$. 
Then Theorem \ref{chain} implies that the $M$-derivative at $a$ of the flexible function $ G $ given by $G(x)= F(x, h(x))$ is well-defined, and satisfies
\begin{equation} \label{chainrule} \frac{d  G}{d_{M} x} (a)\subseteq   \frac{\partial  F}{\partial_{M} x} (a, b)+  \frac{\partial  F}{\partial_{N} y}(a, b)\frac{d  h}{d_{M} x} (a).
\end{equation} 
Because $(a, b)$  is an $(M, N)$-local $L$-minimizer of \eqref{2Lagrange problem}-\eqref{1Lagrange constraint}, it holds that $a$ is an $M$-local $L$-minimizer of $G(x)$. Now $ L $ is stable for $ M  $ and $L\supseteq  {\rm{N} }(F(a, b))$ by assumption, and by \eqref{chainrule} and \eqref{derh} and \eqref{defL} 
\begin{equation}\label{NGL} 
{\rm{N}}\left( 	\frac{d G}{d_{M} x} (a)\right) \subseteq {\rm{N}}\left(   \frac{\partial  F}{\partial_{M} x} (a, b)+  \frac{\partial  F}{\partial_{N} y}(a, b)\frac{d  h}{d_{M} x} (a)\right) \subseteq{\rm{N}}\left(\frac{\partial F}{\partial_{M} x}(a, b)\right) + {\rm{N}}\left( \frac{\gamma_x}{\gamma_y}\frac{\partial F}{\partial_{N} y}(a, b)\right) = L.
\end{equation}
Then by Theorem \ref{necond}	
\begin{equation}\label{inclsion of derivative}
 \frac{d  G}{d_{M} x} (a)\subseteq L.
\end{equation}
Because of \eqref{chainrule} and \eqref{inclsion of derivative} one has 

$$\frac{\partial  F}{\partial_{M} x} (a, b)+  \frac{\partial  F}{\partial_{N} y}(a, b)\frac{d  h}{d_{M} x} (a)\bigcap L \not =\emptyset.$$ 
If two external numbers have non-empty intersection, one of them includes the other, so \eqref{FL} holds directly, or 
\begin{equation}
\label{case 1}  L \subseteq \frac{\partial  F}{\partial_{M} x} (a, b)+  \frac{\partial  F}{\partial_{N} y}(a, b)\frac{d  h}{d_{M} x} (a).
\end{equation}
Then $ \frac{\partial  F}{\partial_{M} x} (a, b)+\frac{\partial  F}{\partial_{N} y}(a, b)\frac{d  h}{d_{M} x} (a) $ is neutricial, so \eqref{NGL} implies that
\begin{equation*}
\frac{\partial  F}{\partial_{M} x} (a, b)+  \frac{\partial  F}{\partial_{N} y}(a, b)\frac{d  h}{d_{M} x} (a)= {\rm{N}}\left( \frac{\partial  F}{\partial_{M} x} (a, b)+  \frac{\partial  F}{\partial_{N} y}(a, b)\frac{\partial  h}{\partial_{M} x} (a) \right) \subseteq L.
\end{equation*}
Hence \eqref{FL} always holds.

Secondly, we prove \eqref{gL}. By \eqref{derh}

$$\frac{\partial g}{\partial x}(a, b)+\frac{\partial g}{\partial y}(a, b)\frac{d h}{d_{M} x}(a)
\subseteq \frac{\partial g}{\partial y}(a, b){\rm{N}}\left(\gamma_x/\gamma_y\right).$$ 
Hence
\begin{equation}\label{ptnp2}  	\lambda	\left( \frac{\partial g}{\partial x}(a, b)+  \frac{\partial g}{\partial y}(a, b)\frac{\partial h}{\partial_{M} x}(a)\right) \subseteq\lambda \frac{\partial g}{\partial y}(a, b){\rm{N}}\left(\gamma_x/\gamma_y\right)
	 \subseteq  \frac{\partial F}{\partial_{N} y}(a, b){\rm{N}}\left(\gamma_x/\gamma_y\right)\subseteq L.
\end{equation} 

Thirdly, we prove \eqref{ptn2}. Observe that $  \frac{\partial  F}{\partial_{N} y}(a) - \lambda \frac{\partial g}{\partial y}(a, b) $ is neutricial by \eqref{ptn1}. Then, using (sub)distributivity
\begin{align*}\label{plusn} \frac{\partial F}{\partial_{M} x}(a, b) - \lambda \frac{\partial g}{\partial x}(a, b)& \subseteq  \frac{\partial  F}{\partial_{M} x}(a, b) - \lambda \frac{\partial g}{\partial x}(a, b)+ \left( \frac{\partial  F}{\partial_{N} y}(a) - \lambda \frac{\partial g}{\partial y}(a, b)\right) \frac{d h}{d_{M} x}(a)\\
	&\subseteq  \frac{\partial  F}{\partial_{M} x}(a, b) +  \frac{\partial  F}{\partial_{N} y}(a)\frac{d h}{d_{M} x}(a)- \lambda \frac{\partial g}{\partial x}(a, b) - \lambda \frac{\partial g}{\partial y}(a, b) \frac{d h}{d_{M} x}(a)\\
	&=\frac{\partial  F}{\partial_{M} x}(a, b) +  \frac{\partial  F}{\partial_{N} y}(a)\frac{d h}{d_{M} x}(a)+\lambda\left(  \frac{\partial g}{\partial x}(a, b)
	+  \frac{\partial g}{\partial y}(a, b) \frac{d h}{d_{M} x}(a)\right). 
\end{align*}

Finally, combining \eqref{ptn2}, \eqref{FL} and \eqref{gL}, we obtain that
\begin{equation*}
\frac{\partial F}{\partial_{M} x}(a, b) - \lambda \frac{\partial g}{\partial x}(a, b)	
\subseteq L+L=L.
\end{equation*}

\end{proof}

\begin{example}Let $ \eps $ be a positive infinitesimal. Consider the optimization problem
	$$F(x, y)=-(1+\eps\oslash)x-(1+\eps\pounds)y^2+\oslash \rightarrow \min$$
	subject to the constraint $x^2+y^2=1$. Let $g(x, y)=1-x^2-y^2$. By direct verification we see that the point $(1/2, \sqrt{3}/2)$ is an $(M, N)$-local $L$-minimizer of the problem, with $M=N=L=\oslash$ and that $g$ and $F$ satisfy at $(1/2, \sqrt{3}/2)$ all remaining assumptions of Theorem \ref{Lagrange multiplier}. Also
	\begin{equation*}	
		\frac{\partial F}{\partial_\oslash x}(1/2, \sqrt{3}/2)=-1+\oslash,\quad  
		\frac{\partial F}{\partial_\oslash y}(1/2, \sqrt{3}/2)=-\sqrt{3}+\oslash, \quad 
		\frac{\partial g}{\partial x}(1/2, \sqrt{3}/2)=1, \quad
		\frac{\partial g}{\partial y}(1/2, \sqrt{3}/2)=\sqrt{3},
	\end{equation*}
and a Lagrange multiplier satisfying the conclusions of Theorem \ref{Lagrange multiplier} is given by $\lambda=1$, for
	\begin{align*} \frac{\partial F}{\partial_\oslash x}(1/2, \sqrt{3}/2)-\lambda \frac{\partial g}{\partial x}(1/2, \sqrt{3}/2)&=\oslash\\
		\frac{\partial F}{\partial_\oslash y}(1/2, \sqrt{3}/2)-\lambda \frac{\partial g}{\partial y}(1/2, \sqrt{3}/2)&=\oslash.
	\end{align*} 
\end{example}


\begin{thebibliography}{99}	
	\bibitem{Van der Corput} van der Corput~JG.  Introduction to the neutrix calculus. Journal d'Analyse  Math\'ematique. 1959;7(1):291-398.
	
	\bibitem{Dienerreeb} Diener~F, Reeb~G. Analyse nonstandard. Hermann; 1989.
	
	\bibitem{Dinis} 	Dinis~B, van den Berg~IP.   Algebraic properties of external numbers. Journal of Logic $\&$ Analysis. 2011;3(9):1-30.
	
	\bibitem{dinisberg 2016} Dinis~B, van den Berg~IP.  Axiomatics for the external numbers of nonstandard analysis. Journal of Logic \& Analysis. 2017;9(7):1-47.
	
	
	\bibitem{DinisBerg} Dinis~B, van den Berg~IP. Neutrices and
	External Numbers. A flexible number system. London: Taylor and Francis; 2019.
	\bibitem{Sequences}  Dinis~B, Tran~VN, van den Berg~IP. {\em On flexible sequences}, Acta Mathematica Vietnamica. 2019;44(4):833-874.
	
	\bibitem{Karl Monty Toda} Karl~JS, Monty~JS, Magdalena~DT. Calculus (sixth edition). Kendall Hunt Publishing Company; 2014. 
	
	\bibitem{KanoveiReeken} Kanovei~V, Reeken~M. Nonstandard analysis, axiomatically. Springer; 2004.
	
	\bibitem{Koudjeti Van den Berg} Koudjeti~F., van den Berg~IP.  Neutrices, external numbers and external calculus. in Nonstandard Analysis in Practice, F. and M. Diener (eds.). Springer Universitext; 145-170, 1995.
	
	
	\bibitem{Lyantsekudryk} Lyantse~W, Kudryk~T.  Introduction to nonstandard analysis. VNTL Publishers, Lviv; 1997.
	
	
	\bibitem{Nelson}  Nelson~E. Internal set theory: a new approach to nonstandard analysis. Bulletin of the American Mathematical Society. 1977;83(6):1165-1198.
	
	\bibitem{Rockafellar} Rockafellar~RT.  Convex analysis. Princeton University Press; 1970.
	
	
	
	\bibitem{Sentivity} Saltelli~A, Ratto~M, Andres~T, Campolongo~F, Cariboni~J,  Gatelli~D, Saisana~M, Tarantola~S. Global sensitivity analysis: The Primer.  Wiley; 2008.
	
	\bibitem{VDBNAA} Van den Berg~IP. Nonstandard Asymptotic Analysis. Springer Lecture Notes in Mathematics 1249; 1987.
	
	\bibitem{Immedecompositionofneutrices} Van den Berg~IP. A decomposition theorem for neutrices. Annals of Pure and Applied Logic. 2010;161(7):851-865.	
	
	
	\bibitem{Taylor} Taylor~JR.  An introduction to error analysis: The study of uncertainties in physical measurements, 2nd ed. University Science Books; 1997.
	
	\bibitem{Nam} Tran~VN, van den Berg~IP.  A parameter method for linear algebra and optimization with uncertainties. Optimization, 2020; 69(1): 21-61.	
\end{thebibliography}
\end{document}